\documentclass[12pt]{article}

\usepackage{amssymb}

\usepackage{amsmath}

\usepackage[mathcal]{euscript}

\usepackage{amsthm}

\usepackage{graphicx}

\usepackage{physics}

\usepackage{float}

\usepackage{tikz}
\usepackage{xcolor}

\newcommand{\HH}{\mathbb{H}}

\newcommand{\XX}{\mathbb{X}}
\newcommand{\YY}{\mathbb{Y}}
\newcommand{\ZZ}{\mathbb{Z}}

\newcommand{\e}{\varepsilon}
\newcommand{\EE}{\mathsf{E}}
\newcommand{\PP}{\mathsf{P}}

\begin{document}

\begin{center}
{\bf \large  Perturbed Markov Chains and \\ Information Networks}
\end{center}

\vspace{1mm}

\begin{center}
Benard Abola, Pitos Seleka Biganda, Sergei Silvestrov, Dmitrii Silvestrov,
Christopher Engstr\"{o}m, John Magero Mango, and
Godwin  Kakuba\footnote{
Benard Abola, M\"{a}lardalen University, V\"{a}ster{\aa}s, Sweden, \\
E-mail: benard.abola@mdh.se \\
Pitos Seleka Biganda, M\"{a}lardalen University, V\"{a}ster{\aa}s, Sweden
and University of Dar es Salaam, Dar es Salaam, Tanzania \\
E-mail: pitos.biganda@mdh.se \\
Sergei Silvestrov, M\"{a}lardalen University, V\"{a}ster{\aa}s, Sweden \\
E-mail: sergei.silvestrov@mdh.se \\
Dmitrii Silvestrov, Stockholm University, Stockholm, Sweden \\
E-mail: silvestrov@math.su.se \\
Christopher Engstr\"{o}m, M\"{a}lardalen University, V\"{a}ster{\aa}s, Sweden \\
E-mail: christopher.engstrom@mdh.se, \\
John Magero Mango,  Makerere University, Kampala, Uganda \\
E-mail: mango@cns.mak.ac.ug \\
Godwin Kakuba, Makerere University, Kampala, Uganda \\
E-mail: godwin.a.kakuba@gmail.com} \\
\end{center}

{\bf Abstract}. The paper is devoted to studies of perturbed Markov chains commonly used for description of information networks. In such models,  the matrix of transition probabilities for the  corresponding Markov chain is usually regularised by adding a special damping matrix multiplied by a small damping (perturbation)  parameter $\e$. We give effective upper bounds for the rate of approximation for stationary distributions of unperturbed Markov chains by stationary distributions of perturbed Markov chains with regularised matrices of transition probabilities, asymptotic expansions for approximating stationary distributions with respect to
damping parameter,  as well as explicit upper bounds for the rate of convergence in ergodic theorems for $n$-step transition probabilities
in triangular array mode, where perturbation parameter $\e \to 0$ and $n \to \infty$,  simultaneously. The results of numerical experiments  are also presented.  \\

{\bf 1. Introduction} \\

Perturbed Markov chains is one of the popular and important objects of studies in the theory of Markov processes and their
applications  to stochastic networks, queuing and reliability models, bio-stochastic systems, and many other stochastic models.

We refer here to some recent books and papers devoted to perturbation problems for Markov type processes,
\cite{AFHo, APZ, BLM, En2, GKP, GS4, KGM, KoLi0, Ko2, MIT, Ni3, NSMD1, Pe2, SiP, SPH,
SiSi1, SiSi2, SiSi3, SiSi4, SiSi5, St9, St10, YZ2, YZ3}. In particular, we would like to mention works   \cite{AFHo, GS4, SiSi1, SiSi2},
where the  extended bibliographies of works in the area and the corresponding methodological and  historical remarks can be found.

We are especially interested in models of Markov chains commonly  used for description of information networks. In such models an information network is represented by the Markov chain associated with  the corresponding node links graph. Stationary distributions and other related characteristics  of information Markov chains usually serve as basic tools for ranking of nodes in information networks.  The ranking problem may be complicated by singularity of the corresponding information Markov chain, where its phase space is split into several  weakly or completely non communicating groups of states. In such models, the  matrix of transition probabilities  $\mathbf{P}_0$ of information Markov chain  is usually regularised  and approximated by  matrix
$\mathbf{P}_{\varepsilon}  = (1 - \varepsilon) \mathbf{P}_0 + \varepsilon \mathbf{D}$, where $\mathbf{D}$ is a so-called damping stochastic matrix with identical rows and all positive elements, while  $\varepsilon \in [0, 1$] is a damping (perturbation) parameter.

The power method  is  often used to approximate the corresponding stationary distribution $\bar{\pi}_{\e}$ by rows of matrix ${\mathbf P}_\e^n$. The damping parameter $0 < \e \leq 1$ should be chosen neither too small nor too large. In the first case, where $\e$  takes too small values, the damping effect will not work against absorbing and pseudo-absorbing effects, since the second eigenvalue for such matrices (determining the rate of convergence in the above mentioned ergodic approximation) take values approaching $1$.   In the second case, the ranking information (accumulated by matrix ${\mathbf P}_0$ via the corresponding stationary distribution) may be partly lost, due to the deviation of  matrix ${\mathbf P}_\e$ from matrix ${\mathbf P}_0$.   This actualises the problem of construction asymptotic expansions for perturbed stationary distribution $\bar{\pi}_{\e}$ with respect to damping parameter $\e$,  as well as studies of asymptotic behaviour of matrices ${\mathbf P}_\e^n$  in triangular array mode, where $\e \to 0$ and $n \to \infty$,  simultaneously.

In this paper, we perform the detailed perturbation analysis of Markov chains with damping component, using the procedure of artificial regeneration for the approximating Markov chains and the coupling method in ergodic theorems for perturbed  regenerative processes. We get effective explicit series representations for the corresponding stationary distributions $\bar{\pi}_\varepsilon$, upper bounds for the deviation $\bar{\pi}_{\varepsilon} - \bar{\pi}_0$, and asymptotic expansions for $\bar{\pi}_{\varepsilon}$ with respect to the perturbation  parameter $\varepsilon$,  as well as get  ergodic theorems and effective  explicit  upper bounds for the rate of convergence in the corresponding ergodic relations for $\mathbf{P}_{\varepsilon}^{n}$, as   $\varepsilon  \to 0$ and $n \to \infty$.  The results of numerical experiments illustrating the above results are also presented.

Real world system consists of interacting units or components. These components constitute what is termed as information networks. With recent advancement in technology, filtering  information has become a challenge in such systems. Moreover, their significance is visible as they find their applications in Internet search engines, biological, financial, transport, queuing networks and many others,  
\cite{Andersson2008},  \cite{AFHo}, \cite{Biganda2017b}, \cite{BrinPage(1998)}, \cite{SS_Eng} -- \cite{SS_EnS4}, \cite{Gleich2015}, \cite{Langville(2011)},  and \cite{SHj2013}. 

PageRank is the link-based criteria that captures the importance of Web pages and provide rankings of the pages in the search engine Google \cite{AFHo},  \cite{Biganda2017b}, \cite{BrinPage(1998)}, \cite{SS_Eng} -- \cite{SS_EnS4}, and \cite{Langville(2011)}. The transition matrix (also called Google matrix $\mathbf{G}$) of a Markov chain in PageRank problem is defined in \cite{Andersson2008}, \cite{Haveliwala2003}, and \cite{Langville(2011)}, as $\mathbf{G} = c\mathbf{P} + (1-c)\mathbf{E}$, where $\mathbf{P}$ is an $n\times n$ row-stochastic matrix (also called hyperlink matrix),
$\mathbf{E}$ (the damping matrix) is the $n\times n$ rank-one stochastic matrix and $c\in(0,1)$ is the damping parameter.

The fundamental concept of PageRank is to use the stationary distribution of the Markov chain on the network to rank Web pages. However, other algorithms similar to PageRank exist in literature, for instance, EigenTrust algorithm \cite{Kamvar2003} and DeptRank algorithm \cite{Battison2012}. In addition, variants of PageRank in relation to some specific networks has been studied, e.g. in \cite{Biganda2018} and  \cite{Biganda2017}; and also updating PageRank due to changes in some network is in literature, for instance, in \cite{Abola2017} and \cite{Abola2018}.

The parameter $c$ is very important in the PageRank definition, because it regulates the level of the uniform noise introduced to the system \cite{AFHo, Langville(2011)}. If $c=1$ there are several absorbing states for the random walk defined by $\mathbf{P}$. However, if $0<c<1$, the Markov chain induced by matrix $\mathbf{G}$ is ergodic \cite{AFHo}. The authors of  \cite{Langville(2011)} argued that the parameter $c$ controls the asymptotic rate of convergence of power method algorithm. Similar arguments were given in  \cite{Andersson2008}, where it is  pointed out that the choice of $c$ is a delicate matter. It may result into convergence and stability problems, if not carefully chosen.

The damping factor $c$ may be denoted and interpreted differently depending on a model being studied. For instance, a model of Markov chain with restart is considered in \cite{APZ}, where  parameter $p$ is the probability to restart the move  and $1-p$ is  the probability to follow the link according to the corresponding transition probability of the above Markov chain. Hence, one can argue that the parameter $p$ has the same interpretation as the damping factor in Google's PageRank problem.

Our representation of perturbed Markov chains is  traditional for perturbed stochastic processes. In fact, PageRank is the stationary distribution of the singularly perturbed Markov chain with perturbation parameter $\e=1-c$. Hence, we wish to point out here that, representation of information network model by a Markov chain with matrix of transition probabilities $\mathbf{P}_{\varepsilon}  = (1 - \varepsilon) \mathbf{P}_0 + \varepsilon \mathbf{D}$ should not create any confusion to the reader.  We perform asymptotic analysis of such Markov chains, in particular,  under the assumption that $\e \to 0$.

The paper  includes  9 sections. In Section 2, we describe the algorithm for  stochastic modelling of  Markov chains with damping component and the procedure of embedding such Markov chains in the model of discrete time regenerative processes with special damping regenerative times. In Section 3, we  present  individual ergodic theorems for Markov chains with damping component and give explicit formulas for the corresponding stationary distributions. In  Section 4, we describe continuity properties of transition probabilities and stationary distributions  with respect to dumping parameter. In Section 5, explicit upper bounds for rates of convergence in approximations of the stationary distributions for Markov chain with damping component are given. In Section 6, we present asymptotic expansions for stationary distribution of  Markov chains with damping component with respect to perturbation (damping) parameter. In Section 7, we describe  coupling algorithms for Markov chains with damping component and get explicit estimates for the rate of convergence in the above mentioned ergodic theorems. In Section 8, we present ergodic theorems for Markov chains with damping component  in a triangular array mode, where time tends to infinity and perturbation parameter tends to zero, simultaneously. In Section 9,   results of numerical experiments, which  show how  results presented in Sections 2--8  can be interpreted and useful in studies of  information networks, are presented, and  some concluding remarks and comments are given. \\

{\bf 2. Markov chains with damping component}  \\

In this section, we introduce the model of Markov chains with damping component, which are used for description of information networks. We also  describe the procedure of embedding such Markov chains in the model of discrete time regenerative processes with special damping regenerative times and present the corresponding renewal type equations for transition probabilities.\vspace{1mm}

{\bf 2.1 Stochastic modelling of Markov chains with damping component}.  Let (a) $\XX = \{1, 2, \ldots, m \}$ be a finite space, (b) $\bar{p} = \langle p_1, \ldots, p_m \rangle$, $\bar{d} = \langle d_1, \ldots, d_m \rangle$, and $\bar{q} = \langle q_0, q_1 \rangle$ be three discrete probability distributions, (c) $\mathbf{P}_0 = \| p_{0, ij} \|$ be a  $m \times m$ stochastic matrix and  $\mathbf{D} = \| d_{ij} \|$ be a  $m \times m$ damping stochastic matrix with elements $d_{ij} = d_j > 0, i, j = 1, \ldots,  m$, and $\mathbf{P}_\e  = \| p_{\e, ij} \| = (1- \e)\mathbf{P}_0 + \e\mathbf{D}$ is a stochastic matrix with elements $p_{\e, ij} = (1- \e) p_{0, ij} + \e d_j, i, j = 1, \ldots, m$, where $\e \in [0, 1]$.

 Let us first describe an algorithm for stochastic modelling of a discrete time,  homogeneous Markov chain $X_{\e, n}, n = 0, 1, \ldots$, with the phase space $\XX$,
 the initial distribution $\bar{p}$, and the matrix of transition probabilities $\mathbf{P}_\e$.

Let (a) $U$ be a random variable taking values in space $\XX$ and such that $\PP \{U = j  \} = p_j, j \in \XX$; (b) $U_{i, n}, i \in \XX, n = 1, 2,  \ldots$ be a family of independent random variables taking values in space $\XX$ and such that $\PP \{U_{i, n} = j  \} = p_{0, ij}, i, j \in \XX, n = 1, 2, \ldots$;  (c) $V_n, n = 0, 1, \ldots$ be a sequence of independent random variables taking values in space $\XX$ and such that $\PP \{V_n = j  \} = d_j, j \in \XX,  n = 1, 2, \ldots$   (d) $W$ is a binary random variable taking values $0$ and $1$ with probabilities, respectively $q_0$ and $q_1$; (e) $W_{\e, n}, n = 1, 2,  \ldots$ be, for every $\e \in [0, 1]$,  a sequence of  independent binary random variables taking values $0$ and $1$ with probabilities, respectively, $1 - \e$ and $\e$, for $n = 1, 2, \ldots$; (d) the random variables $U, W$, the
family of random variables  $U_{i, n}, i \in \XX, n = 1, 2, \ldots$, and the random sequences $V_{\e, n}, n = 1, 2, \ldots$ and $W_{\e, n}, n = 1,, 2, \ldots$ are mutually independent, for every $\e \in [0, 1]$.

Let us now define, for every $\e \in [0, 1]$, the random sequence $X_{\e, n}, n = 0, 1, \ldots$, by the following recurrent relation,
\begin{equation}\label{ertea}
X_{\e, n} = U_{X_{\e, n -1}, n} {\rm I}(W_{\e, n} = 0) + V_{n} {\rm I}(W_{\e, n} = 1),  n = 1, 2, \ldots, \,  X_{\e, 0} = U.
\end{equation}

It is readily seen that the random sequence $X_{\e, n}, n = 0, 1, \ldots$ is, for every $\e \in [0, 1]$, a homogeneous Markov chain with phase space $\XX$, the initial distribution $\bar{p}$ and the matrix of transition probabilities $\mathbf{P}_\e$. \vspace{1mm}

{\bf 2.2 Regenerative properties of Markov chains with damping component}. Let us consider the extended random sequence,
\begin{equation}\label{erteaba}
Y_{\e,n} = (X_{\e, n}, W_{\e, n}), n = 1, 2, \ldots, \, X_{\e, 0} = U, W_{\e, 0} = W.
\end{equation}

This random sequence also is,  for every $\e \in [0, 1]$, a
homogeneous Markov chain,  with phase space $\YY = \XX \times \{0, 1 \}$, the initial distribution $\overline{pq} = \langle p_i q_r$, $(i, r) \in \YY \rangle$ and  transition probabilities,
 \begin{equation*}
 p_{\e, ir, jk} = \PP \{ X_{\e, 1} = j, W_{\e, 1} = k / X_{\e, 0} = i, W_{\e, 0} = r \} \makebox[5mm]{}
 \end{equation*}
 \begin{equation}\label{gotyr}
 = \left\{
 \begin{array}{lll}
 (1-  \e) p_{0, ij}  & \text{for} \  i, j \in \XX, \,  r = 0, 1, \,  k = 0, \\
 \e d_j   & \text{for} \  i, j \in \XX, \,  r = 0, 1, \,  k = 1.
 \end{array}
 \right.
 \end{equation}

 It is worth  to note that the transition probabilities $p_{\e, ir, jk} =  p_{\e, i, jk}, \, (i, r)$, $(j, k) \in \YY$ do not depend on $r \in \XX$ and on $i \in \XX$ if $k = 1$.

 Let us, assume that the dumping (perturbation) parameter $\e \in (0, 1]$.

 Let us define times of sequential hitting  state $1$
 by the second component $W_{\e, n}$,
 \begin{equation}\label{gotyrta}
T_{\e, n} = \min(n > T_{\e, n-1},  W_{\e, n} = 1), n = 1, 2, \ldots, \, T_{\e, 0} = 0.
  \end{equation}

The random sequence  $Y_{\e, n}, n = 0, 1, \ldots$ is also a discrete time regenerative process with ``damping''  regeneration
 times $T_{\e, n} , n = 0, 1, \ldots$.

It follows from independence of transition probabilities $p_{\e, ir, jk}$ of the Mar\-kov chain $Y_{\e,n}$, given by relation (\ref{gotyr}), on 
$(i, r) \in \YY$ if $k = 1$.

 This is a standard regenerative process, if the initial distribution  $\overline{pq} = \overline{dq}_1  = \langle d_i q_{1, r}$, $(i, r)  \in \YY \rangle$, where $\bar{q}_1 = \langle q_{1, 0} = 0, q_{1,1} = 1 \rangle$.

Otherwise,  $Y_{\e, n}$ is a regenerative process with the transition period $[0, T_{\e, 1})$.

It is also useful to note that the inter-regeneration times $S_{\e, n} = T_{\e, n} - T_{\e, n-1}, n = 1, 2, \ldots$ are i.i.d. geometrically distributed random variables, with
parameter $\e$, i.e.,
\begin{equation}\label{geom}
\PP \{ S_{\e, 1} = n \} = \left\{
\begin{array}{lll}
 0 & \text{for} \  n = 0, \\
 \e (1 - \e)^{n-1} & \text{for} \  n = 1, 2, \ldots.
 \end{array}
 \right.
 \end{equation}

 {\bf 2.3 Renewal type equations for transition probabilities of the  Markov chains with damping component}.  Let us denote by ${\cal P}_m$ the class of all initial distributions   $\bar{p} = \langle p_i , i \in \XX \rangle$.

 Let us introduce $n$-step transition probabilities
 for the Markov chain $X_{\e, n}$, for $ i, j \in \XX, \, n = 0, 1, \ldots$,
\begin{equation}\label{trans}
p_{\e, ij}(n) = \PP \{ X_{\e, n} = j / X_{\e, 0} = i \},
\end{equation}
and $n$-sep probabilities, for $\bar{p} \in {\cal P}_m,  j \in \XX, \, n = 0, 1, \ldots$,
\begin{equation}\label{transa}
p_{\e, \bar{p}, j}(n) =  \PP_{\bar{p}} \{ X_{\e, n} = j \} = \sum_{i \in \XX} p_i p_{\e, ij}(n).
\end{equation}

Here and henceforth, symbols $\PP_{\bar{p}}$ and $\EE_{\bar{p}}$ are used for probabilities and expectations  related to a Markov chain with an initial distribution
$\bar{p}$. In the case, where the initial distribution is concentrated in a state $i$ the above symbols take the forms $\PP_i$ and $\EE_i$.

Obviously, $p_{\e, ij}(0) = {\rm I} (i = j), i, j \in \XX$ and $p_{\e, ij}(1) = p_{\e, ij}, i, j \in \XX$. Also, $p_{\e, \bar{p},  j}(0) = p_j,  j \in \XX$.

Let us also denote by ${\cal PQ}_m$ the class of all initial distributions   $\overline{pq} = \langle p_iq_r , (i, r) \in \YY \rangle$.

Analogously, let us introduce $n$-step transition probabilities for the Markov chain $Y_{\e, n}$, for $ \, (i, r), (j, k) \in \YY, \, n = 0, 1, \ldots$,
\begin{equation}\label{transbaw}
p_{\e, ir, jk}(n) = \PP \{ Y_{\e, n} = (j, k) / Y_{\e, 0} = (i, r) \},
\end{equation}
and $n$-sep probabilities, for $\overline{pq} \in {\cal PQ}_m,  (j, k) \in \YY, \, n = 0, 1, \ldots$,
\begin{equation}\label{transabaw}
p_{\e, \overline{pq}, jk}(n) =  \PP_{\overline{pq}} \{ Y_{\e, n} = (j, k) \} = \sum_{(i, r) \in \YY} p_iq_r p_{\e, ir, jk}(n)
\end{equation}

Obviously, $p_{\e, ir, jk}(0)  = {\rm I} ((i, r) = (j, k)), (i, r), (j, k) \in \YY$ and $p_{\e, ir, jk}(1)  = p_{\e, ir, jk}, (i, r), (j, k) \in \YY$. Also, $p_{\e, \overline{pq},  (j, k)}(0) = p_jq_k,  (j, k) \in \YY$.

Independence of  the transition probabilities $p_{\e, ir, jk} =  p_{\e, i, jk}, \, (i, r), (j, k) \in \YY$ on $r \in \XX$ and on $i \in \XX$  if $k = 1$, implies  that $n$-step
transition probabilities $p_{\e, ir, jk}(n) = p_{\e, i, jk}(n), (i, r), (j, k) \in \YY, n = 0, 1, \ldots$ also are independent on $r \in \XX$  and on $i \in \XX$ if $k = 1$.

This also implies that  $n$-step probabilities $p_{\e, \overline{pq}, jk}(n) = p_{\e, \overline{p}, jk}(n), \overline{pq} \in {\cal PQ}_m, (j, k) \in \YY, n = 1, 2, \ldots$ are independent on the initial distribution $\bar{q}$.

Let us  assume that the initial distribution  $\overline{pq} = \overline{dq}_1$.  As was mentioned above, $Y_{\e, n}$ is, in this case,  the  standard regenerative process with regeneration times $T_{\e, n}, n = 0, 1, \ldots$.

This fact and relations (\ref{gotyr})  and  (\ref{geom}) imply that
probabilities  $p_{\e, \overline{pq}_1, jk}(n), n = 0, 1, \ldots$ are, for every $j \in \XX, \, k = 0, 1$,  the unique bounded solution for the following discrete time renewal equation,
 \begin{equation}\label{gotyrtbaramo}
p_{\e, \overline{dq}_1, jk}(n) = q_{\e,  \overline{dq}_1, ik}(n) + \sum_{l = 1}^n  p_{\e,  \overline{dq}_1, jk}(n - l) \e (1- \e)^{l- 1}, n \geq 0,
\end{equation}
where, for  $j \in \XX, k = 0, 1, n \geq 0$,
\begin{equation*}
q_{\e,  \overline{dq}_1, jk}(n) = \PP_{ \overline{dq}_1} \{Y_{\e, n} = (j, k), T_{\e, 1} > n \} \makebox[40mm]{}
\end{equation*}
\begin{equation}\label{qera}
= \left\{
\begin{array}{lll}
p_{0, \bar{d}, j}(n)(1- \e)^n {\rm I}(n > 0) & \text{if} \  k = 0, \makebox[1mm]{} \vspace{1mm} \\
d_j {\rm I}(n = 0) & \text{if} \  k = 1.
\end{array}
\right.
\end{equation}

Let us  now consider the general case, with some initial distribution  $\overline{pq} = \langle p_i q_r, (i, r)  \in \YY \rangle \in {\cal PQ}_m$.
As was mentioned above, $Y_{\e, n},$ is, in this case,  the regenerative process with regeneration times $T_{\e, n}, n = 0, 1, \ldots$ and the transition period
$[0, T_{\e, 1})$.

 In this case, probabilities $p_{\e,  \overline{pq}, jk}(n)$ and
$p_{\e,  \overline{dq}_1, jk}(n)$ are, for $j \in \XX, k = 0, 1$, connected by the following renewal type relation,
\begin{equation}\label{gotyrtbarakabav}
p_{\e,  \overline{pq}, jk}(n) = q_{\e, \overline{pq}, jk}(n)   + \sum_{l = 1}^n  p_{\e,  \overline{dq}_1, jk}(n - l) \e (1- \e)^{l- 1}, \, n \geq 0,
\end{equation}
where, for  $j \in \XX, k = 0, 1, n \geq 0$,
\begin{equation*}
q_{\e,  \overline{pq}, jk}(n) = \PP_{ \overline{pq}} \{Y_{\e, n} = (j, k), T_{\e, 1} > n \} \makebox[40mm]{}
\end{equation*}
\begin{equation}\label{qeran}
= \left\{
\begin{array}{lll}
p_{0, \bar{p}, j}(n)(1- \e)^n {\rm I}(n > 0) & \text{if} \  k = 0, \makebox[1mm]{} \vspace{1mm}\\
p_j {\rm I}(n = 0) & \text{if} \  k = 1.
\end{array}
\right.
\end{equation}

The summation of renewal equations (\ref{gotyrtbaramo}) over $k = 0, 1$  yields  the discrete time renewal equation for  probabilities  $p_{\e, \bar{d}, j}(n), n = 0, 1, \ldots$, which are, for every $j \in \XX$, the unique bounded solution for this  equation,
\begin{equation}\label{gotyrtbaraka}
p_{\e, \bar{d}, j}(n) = p_{0, \bar{d}, j}(n)(1- \e)^n   + \sum_{l = 1}^n  p_{\e, \bar{d}, j}(n - l) \e (1- \e)^{l - 1}, \, n \geq 0.
\end{equation}

Also, the summation of renewal type equations (\ref{gotyrtbarakabav}) over $k = 0, 1$ yields that,
in the case of general initial distribution  $\overline{p} = \langle p_i, i  \in \XX \rangle$, the probabilities $p_{\e, \bar{p}, j}(n)$ and $p_{\e, \bar{d}, j}(n)$
are, for every $j \in \XX$, connected by the following renewal type relation,
\begin{equation}\label{gotyrtbarakaba}
p_{\e, \bar{p}, j}(n) = p_{0, \bar{p}, j}(n)(1- \e)^n   + \sum_{l = 1}^n  p_{\e, \bar{d}, j}(n - l) \e (1- \e)^{l - 1}, \, n \geq 0.
\end{equation}

{\bf 3. Stationary distributions of Markov chains with \\
\makebox[10mm]{}  damping component} \\

In this section we present ergodic relations for transition probabilities of Markov chains with damping component. \vspace{1mm}

{\bf 3.1 Stationary distributions of the  Markov chains $X_{\e, n}$ and $Y_{\e, n}$}.
Let us  describe ergodic properties of Markov chains $X_{\e, n}$ and $Y_{\e, n}$, for the case, where $\e \in (0, 1]$. \vspace{1mm}

{\bf Lemma 1}. {\em Let $\e \in (0, 1]$. Then the following ergodic relation takes place for any initial distribution $\overline{pq} \in {\cal PQ}_m$ and $(j, k) \in \YY$,}
 \begin{equation}\label{gotyrtbaramoto}
p_{\e, \overline{pq}, jk}(n) \to \pi_{\e, jk} = \e \sum_{l = 0}^{\infty}  q_{\e, \overline{dq}_1, jk}(l)  \ {\rm as } \ n \to \infty.
 \end{equation}

{\bf Proof}. The geometrical distribution of the inter-regeneration time $S_{\e, 1}$ is aperiodic and has expectation $\e^{-1}$.

This makes it possible  to apply the discrete time
renewal theorem (see, for example, \cite{Fel}) to the renewal equation (\ref{gotyrtbaramo}) that  yields the following ergodic relation, for  $(j, k) \in \YY$,
 \begin{equation}\label{gotyo}
p_{\e, \overline{dq}_1, jk}(n) \to \pi_{\e, jk} \ {\rm as} \ n \to \infty.
 \end{equation}

 Obviously $q_{\e, \overline{pq}, jk}(n)  \to 0$ as $n \to \infty$, for $(j, k) \in \YY$.

 Let us also define $p_{\e, \overline{dq}_1, jk}(n - l) = 0$, for $l > n$. Relation (\ref{gotyo}) implies that $p_{\e, \overline{dq}_1, jk}(n - l) \to \pi_{\e, jk}$ as $n \to \infty$,
 for $l \geq 0$ and $(j, k) \in \YY$.

Using relation (\ref{gotyrtbarakabav}), the latter two asymptotic relations and the Lebesgue theorem, we get, for  $\overline{pq} \in {\cal PQ}_m,  (j, k) \in \YY$,
\begin{align}\label{gotyrtbarakabavo}
\lim_{n \to \infty} p_{\e, \overline{pq}, jk}(n) & = \lim_{n \to \infty}  q_{\e, \overline{pq}, jk}(n)  \vspace{2mm} \nonumber \\
&  + \lim_{n \to \infty} \sum_{l = 1}^\infty  p_{\e, \overline{dq}_1, jk}(n - l) \e (1- \e)^{l - 1} = \pi_{\e, jk}.
\end{align}

The proof is complete. $\Box$

The following lemma is the direct corollary of Lemma 1. \vspace{1mm}

{\bf Lemma 2}. {\em Let $\e \in (0, 1]$. Then the following ergodic relation takes place for any initial distribution $\bar{p} \in {\cal P}_m$ and $j \in \XX$,}
 \begin{equation}\label{gotyrtbaramotok}
p_{\e, \bar{p}, j}(n) \to \pi_{\e, j} = \e \sum_{l = 0}^{\infty}  p_{0, \bar{d}, j}(l)(1- \e)^l \ {\rm as } \ n \to \infty.
 \end{equation}

 It is useful to note that, for every $\e \in (0, 1]$,  the phase space $\XX$ is one communicative, aperiodic class of  states for the Markov chain $X_{\e, n}$, and its stationary distribution $\bar{\pi}_\e = \langle \pi_{\e, j} , j \in \XX \rangle$ is the unique positive solution for the system of linear equations,
\begin{equation}\label{ergom}
\sum_{i \in \XX} \pi_{\e, i} p_{\e, ij} = \pi_{\e, j}, j \in \XX, \ \sum_{j \in \XX} \pi_{\e, j} = 1.
\end{equation}

Also, the stationary probabilities $\pi_{\e, j}$ can be represented in the form $\pi_{\e, j}  = e_{\e, j}^{-1}, ,j \in  \XX$, via the expected return times $e_{\e, j}$, with the use of regeneration property of the Markov 
chain  $X_{\e, n}$ at moments of return in  state $j$.

The series representation (\ref{gotyrtbaramoto}) for the stationary distribution of Markov chain $X_{\e, n}$ is based on the use of alternative damping regeneration times. This representation is, by our opinion,  a more effective 
tool for performing asymptotic perturbation analysis for Markov chains with damping component.

{\bf 3.2 The stationary distribution of the Markov chain $X_{0, n}$}. Let us  describe ergodic properties of Markov chains $X_{0, n}$. Its
ergodic properties  are determined by communicative properties its phase space $\XX$ and the matrix of
transition probabilities $\mathbf{P}_0$. The simplest case is where the following condition holds:
\begin{itemize}
\item [${\bf A}_1$:] The phase space $\XX$ is one aperiodic class of  communicative states for the Markov chain $X_{0, n}$.
\end{itemize}

In this case, the following ergodic relation holds, for  any $\bar{p} \in {\cal P}_m, j \in \XX$,
\begin{equation}\label{erg}
p_{0, \bar{p}, j}(n)  \to \pi_{0, j} \ {\rm as} \ \e \to 0,
\end{equation}

The stationary distribution $\bar{\pi}_0 = \langle \pi_{0, j}, j \in \XX \rangle$ is the unique positive solution of the system for linear equations,
\begin{equation}\label{ergsys}
\sum_{i \in \XX} \pi_{0, i} p_{0, ij} = \pi_{0, j}, j \in \XX, \ \sum_{j \in \XX} \pi_{0, j} = 1.
\end{equation}

A more complex is the case, where the following condition holds, for some $h > 1$:
\begin{itemize}
\item [${\bf B}_1$:] The phase space $\XX = \cup_{j = 1}^h \XX^{(j)}$, where $X^{(j)}, j = 1, \ldots, h$ are non-intersecting subsets of $\XX$, (b) $X^{(j)}, j = 1, \ldots, h$ are non-empty, closed aperiodic  classes of  communicative  states for the Markov chain $X_{0, n}$.
\end{itemize}

If the initial distribution of the Markov chain $X_{0, n}$ is concentrated at the set $\XX^{(j)}$, for some $j = 1, \ldots, h$, then $X_{0, n} = X^{(j)}_{0, n}, n = 0, 1, \ldots$ can be considered as a Markov chain with the reduced phase space
$\XX^{(j)}$ and  the matrix of transition probabilities $\mathbf{P}_{0, j} = \|�p_{0, rk} \|_{k, r \in \XX^{(j)}}$.

According to condition ${\bf B}_1$, there exists, for any $r, k \in \XX^{(j)}, j = 1, \ldots, h$,
\begin{equation}\label{starew}
p_{0, r k}(n) \to \pi^{(j)}_{0, k} \ {\rm as} \ n \to \infty,
\end{equation}
where $\bar{\pi}_0^{(j)} = \langle \pi^{(j)}_{0, k}, k \in \XX^{(j)} \rangle$ is, for $j = 1, \ldots, h$, the stationary distribution of the Markov chain $X^{(j)}_{0, n}$.

The stationary distribution  $\bar{\pi}_0^{(j)}$ is, for every $j = 1, \ldots, h$, the unique positive solution for the system of linear equations,
\begin{equation}\label{statoifabe}
\pi^{(j)}_{0, k}  = \sum_{r \in \XX^{(j)}} \pi^{(j)}_{0, r} p_{0, rk}, k \in \XX^{(j)}, \ \sum_{k \in \XX^{(j)}}  \pi^{(j)}_{0, k} = 1.
\end{equation}

Let us denote, for $\bar{p} \in {\cal P}_m, j = 1, \ldots, h$,
\begin{align}\label{recoupta}
f^{(j)}_{\bar{p}}  & = \sum_{i \in \XX^{(j)}} p_i .
\end{align}
\vspace{1mm}

{\bf Lemma 3}. {\em Let condition ${\bf B}_1$ holds. Then, the following ergodic relation takes place, for $\bar{p} \in {\cal P}_m$ and
$k \in \XX^{(j)}, j = 1, \ldots, h$,}
\begin{equation}\label{convert}
p_{0, \bar{p}, k}(n) \to  \pi_{0, \bar{p},  k}  =
f^{(j)}_{\bar{p}}   \pi^{(j)}_{0, k}  \ {\rm as} \ n \to \infty.
\end{equation}

\vspace{1mm}

{\bf Remark 1}.  Ergodic relation {\rm (\ref{convert})} shows that in the singular case, where condition ${\bf B}_{1}$ holds,
the stationary probabilities $\pi_{0, \bar{p},  k}$ defined by the asymptotic relation {\rm  (\ref{convert})} may depend on the initial distribution.  
The stationary distributions $\bar{\pi}_{0, \bar{p}} = \langle \pi_{0, \bar{p},  k}, k \in \XX \rangle$ and  $\bar{\pi}_{0, \bar{d}} = 
\langle \pi_{0, \bar{d},  k}, k \in \XX \rangle$ coincide, if probabilities $f^{(g)}_{\bar{p}}  = f^{(g)}_{\bar{d}}, \, g = 1, \ldots, h$.  These relations obviously hold for any initial distribution $\bar{p} \in {\cal P}_m$ in the regular case, where condition ${\bf A}_{1}$ holds.  \\

{\bf 4. A perturbation model for Markov chains with  \\
\makebox[11mm]{} damping component}. \\

In this section, we show in which way Markov chains with damping component can be interpreted as a stochastic perturbed model. We also present results concerned with continuity of  stationary distributions $\bar{\pi}_\e$ with respect to damping (perturbation) parameter $\e \to 0$.
\vspace{1mm}

{\bf 4.1 Continuity property transition probabilities}.  In what follows,  relation $\e \to 0$ is a reduced version of relation $0 < \e \to 0$.

The Markov chain $X_{\e, n}$ has the matrix of transition probabilities  $\mathbf{P}_\e = (1-\e)\mathbf{P}_0 + \e \mathbf{D}$. Obviously, for $i, j \in \XX$,
\begin{equation}\label{obvio}
p_{\e, ij} \to p_{0, ij} \ {\rm  as} \ \e \to 0.
\end{equation}

Also, as well known, matrix $\| p_{\e, ij}(n) \| =  \mathbf{P}_\e^n$, for $n = 0, 1, \ldots$, where $\mathbf{P}_\e^0 = \| {\rm I}(i = j) \|$.

Therefore, the following asymptotic relation holds, for $n \geq 0, i, j \in \XX$,
\begin{equation}\label{transava}
p_{\e, ij}(n) \to p_{0, ij}(n) \ {\rm as} \ \e \to 0.
\end{equation}

This relation let one consider the Markov chain $X_{\e, n}$, for $\e \in (0, 1]$,   as a perturbed version of  the
Markov chain $X_{0, n}$ and to interpret the damping parameter $\e$ as a perturbation parameter.

It was mentioned above,  that the phase space $\XX$ of the perturbed Markov chain  $X_{\e, n}$ is one communicative, aperiodic class of states, for every $\e \in (0, 1]$.

As far as the unperturbed Markov chain $X_{0, n}$ is concerned, there are two different cases.

The first one is, where condition ${\bf A}_1$ holds, i.e., the phase space $\XX$ is also  one communicative class of states for the Markov chain $X_{0, n}$. In this case, one can refer to
the corresponding perturbation model as regular.

The second one is, where condition ${\bf B}_1$ holds, i.e., the phase space $\XX$ is not  one communicative class of states for the Markov chain $X_{0, n}$. In this case, one can refer to
the corresponding perturbation model as singular.
\vspace{1mm}

{\bf 4.2 Continuity property of  stationary distributions for regularly perturbed Markov chains with damping component}. The following proposition takes place. \vspace{1mm}

{\bf Lemma 4}. {\em Let condition ${\bf A}_1$  holds. Then, the following asymptotic relation holds, for $j \in \XX$,}
\begin{equation}\label{reprenasas}
\pi_{\e, j}  \to \pi_{0, j}  \ {\rm as} \ \e \to 0.
\end{equation}

{\bf Proof}. Let $\nu_\e$ be a  random variable geometrically distributed with parameter $\e$, i.e., $\PP \{ \nu_\e = n \} = \e (1 - \e)^{n-1}, n = 1, 2, \ldots$. Obviously,
\begin{equation}\label{reprenast}
\nu_\e - 1 \stackrel{\PP}{\longrightarrow} \infty \ {\rm as}  \  \e \to 0.
\end{equation}

In the case, where condition ${\bf A}_1$ holds, we get using  relations (\ref{erg})  and (\ref{reprenast})  that the following asymptotic
relation holds, for $j \in \XX$,
\begin{equation}\label{reprebasan}
p_{0, \bar{d}, j}(\nu_\e -1) \stackrel{\PP}{\longrightarrow} \pi_{0, j} \ {\rm as} \ \e \to 0.
\end{equation}

It is readily seen that the following representation takes place for the stationary probabilities $\pi_{\e, j}, i \in \XX$,
\begin{equation}\label{repren}
\pi_{\e, j} = \e \sum_{l = 0}^{\infty}  p_{0, \bar{d}, j}(l)(1- \e)^l = \EE p_{0, \bar{d}, j}(\nu_\e -1).
\end{equation}

Since the sequence  $p_{\bar{p}, j}(n), n = 0, 1, \ldots$ is a bounded, relations (\ref{reprebasan}),  (\ref{repren}) and the corresponding variant of the Lebesgue theorem imply that the
following asymptotic relation holds, for $j \in \XX$,
\begin{equation}\label{reprenas}
\pi_{\e, j} =  \EE p_{0, \bar{d}, j}(\nu_\e -1) \to \pi_{0, j} \ {\rm as} \ \e \to 0.
\end{equation}

The proof is complete. $\Box$ \vspace{1mm}

{\bf 4.3 Continuity property of  stationary distributions for singularly perturbed Markov chains with damping component}.  In this case the following lemma takes place. \vspace{1mm}

{\bf Lemma 5}. {\em Let condition  ${\bf B}_1$ holds. Then, the following asymptotic relation holds, for $k \in \XX$,}
\begin{equation}\label{reprenasasa}
\pi_{\e, k}  \to \pi_{0, \bar{d}, k}  \ {\rm as} \ \e \to 0.
\end{equation}

{\bf Proof}. In the case, where condition ${\bf B}_1$ holds, we get using  relations (\ref{convert})  and (\ref{reprenast})  that the following asymptotic relation holds,
\begin{equation}\label{reprebasabaw}
p_{0, \bar{d}, k}(\nu_\e -1) \stackrel{\PP}{\longrightarrow} \pi_{0,\bar{d},  k} \ {\rm as} \ \e \to 0, \ {\rm for} \  k \in \XX.
\end{equation}

Analogously to the relation (\ref{reprenas}), we get, using  relatins (\ref{reprebasan})  and (\ref{reprebasabaw}), the
following asymptotic relation, for $k \in \XX$,
\begin{equation}\label{reprenasnao}
\pi_{\e, k} =  \EE p_{0, \bar{d}, k}(\nu_\e -1) \to \pi_{0, \bar{d}, k} \ {\rm as} \ \e \to 0.
\end{equation}

The proof is complete. $\Box$ \vspace{1mm}

{\bf Remark 2}. Lemmas {\rm 4} and  {\rm 5}  imply that, in the case where condition ${\bf A}_1$ holds, the continuity property for stationary distributions $\bar{\pi}_\e$  {\rm (}as $\e \to 0${\rm )} takes place. However, in the case where condition ${\bf B}_1$ holds, the continuity property for stationary distributions $\bar{\pi}_\e$  {\rm (}as $\e \to 0${\rm )} takes place under the
additional assumption under  that  $f^{(g)}_{\bar{p}}  = f^{(g)}_{\bar{d}}, g = 1, \ldots, h$.   \\

{\bf 5. Rate of convergence for perturbed stationary distributions} \\

In this section, we obtain explicit upper bounds for deviations of stationary distributions of Markov chains $X_{\e, n}$  and $X_{0, n}$.

\vspace{1mm}

{\bf 5.1 Rate of convergence for stationary distributions of regularly perturbed Markov chains with damping component}.
Let us get some explicit upper bound for the rate of convergence in the asymptotic relation (\ref{reprenasas}) for the case, where  condition ${\bf A}_1$ holds.

It is well known that, under condition ${\bf A}_1$,  the rate of convergence in the ergodic relation (\ref{reprenasas}) is exponential. This means that
there exist some constants $C = C(\mathbf{P}_0) \in [0, \infty)$, $\lambda = \lambda(\mathbf{P}_0) \in [0, 1)$, and distribution
$\bar{\pi}_0 = \langle \pi_{0, j}, j \in \XX \rangle$, with all positive components such that the following relation holds,
 \begin{equation}\label{asnao}
 \max_{i,j \in \XX} |p_{0, ij}(n) - \pi_{0, j} | \leq C \lambda^n, \ n \geq 1.
\end{equation}

In fact, condition ${\bf A}_1$ is equivalent to the following condition:
\begin{itemize}
\item [${\bf A}_2$:] There exist a constants $C = C(\mathbf{P}_0) \in [0, \infty)$, $\lambda = \lambda(\mathbf{P}_0) \in [0, 1)$, and a distribution
$\bar{\pi}_0 = \langle \pi_{0, j}, j \in \XX \rangle$ with all positive components such that relation (\ref{asnao}) holds.
\end{itemize}

Indeed, condition ${\bf A}_2$ implies that probabilities $p_{0, ij}(n) > 0, i, j \in \XX$ for all large enough $n$. This implies that $\XX$ is
one  aperiodic class of communicative  states.  Also, condition ${\bf A}_2$ implies that $p_{0, ij}(n)  \to  \pi_{0, j}$ as $n \to \infty$, for $i, j \in \XX$, and, thus,
$\bar{\pi}_0 $ is the stationary distribution for the Markov chain $X_{0, n}$.

According to the Perron-Frobenius theorem, the role of $\lambda$ can play  the absolute value of  the second (by absolute value),  eigenvalue for matrix $\mathbf{P}_0$. As far as constant $C$ is concerned, we refer to the book \cite{Fel}, where one can find the algorithms
which let one compute this constant.

An alternative and more simple variant of inequalities appearing in condition  ${\bf A}_2$ may be obtained with the use of coupling method. We give them in Section 8.

The following theorem present explicit upper bounds for deviation of stationary distributions for
Markov chains $X_{\e, n}$  and $X_{0, n}$. \vspace{1mm}

{\bf Theorem 1}. {\em Let condition ${\bf A}_2$ holds. Then the following relation holds, for $j \in \XX$, }
\begin{equation}\label{estimak}
| \pi_{\e, j}  - \pi_{0, j} | \leq  \e ( | d_j - \pi_{0, j}|  + \frac{C \lambda}{1 - \lambda}).
\end{equation}

{\bf Proof}. The inequalities appearing in condition  ${\bf A}_2$ imply that the following relation holds, for $n \geq 1, j \in \XX$,
\begin{align}\label{estimaba}
|p_{0, \bar{d}, j}(n) - \pi_{0, j} |  & = | \sum_{i \in \XX} (d_i p_{0, ij}(n) - d_i \pi_{0, j}) |  \vspace{2mm} \nonumber \\
& \leq \sum_{i \in \XX} d_i | p_{0, ij}(n)  -  \pi_{0, j})| \leq C \lambda^n.
\end{align}

Using relations (\ref{gotyrtbaramoto}) and (\ref{estimaba}), we get the following
estimate,  for $j \in \XX$,
\begin{align}\label{esti}
| \pi_{\e, j}  - \pi_{0, j} | & \leq | \e \sum_{l = 0}^{\infty}  p_{0, \bar{d}, j}(l)(1- \e)^l -  \pi_{0, j}  | \vspace{2mm} \nonumber \\
& = | \e \sum_{l = 0}^{\infty}  p_{0, \bar{d}, j}(l)(1- \e)^l - \e \sum_{l = 0}^{\infty}  \pi_{0, j} (1- \e)^l | \vspace{2mm} \nonumber \\
& \leq \e | d_j - \pi_{0, j}| + \e \sum_{l = 1}^{\infty} C\lambda^l (1- \e)^l  \vspace{2mm} \nonumber \\
& \leq  \e ( | d_j - \pi_{0, j}| +  \frac{C \lambda (1-\e)}{1 - \lambda (1 - \e)})  \vspace{2mm} \nonumber \\
& \leq  \e ( | d_j - \pi_{0, j}|  + \frac{C \lambda}{1 - \lambda}).
\end{align}

The proof is complete. $\Box$ \vspace{1mm}

{\bf Remark 3}. The quantities $| d_j - \pi_{0, j}|$ appearing in inequality {\rm (\ref{estimak})} is, in some sense, determined by  a prior information about the stationary probabilities $\pi_{0, j}$. They take smaller values if one can choose initial distribution $\bar{p}$ with smaller deviation from the stationary distribution $\bar{\pi}_0$. Inequalities $| d_j - \pi_{0, j}| \leq d_j \vee (1 - d_j) \leq 1$ let one replace the term   $| d_j - \pi_{0, j}|$ in inequality {\rm (\ref{estimak})} by quantities independent on the corresponding stationary probabilities $\pi_{0, j}$. \vspace{1mm}

{\bf Remark 4}. Theorem {\rm 1} remains valid even we weaken condition ${\bf A}_2$ by omitting in it the assumption of positivity for the distribution  $\bar{\pi}_0 = \langle \pi_{0, i}, i \in \XX \rangle$. In this case, condition ${\bf A}_2$ implies that the phase space $\XX = \XX_1 \cup \XX_0$, where $\XX_1 = \{ i \in \XX: \pi_{0, i} > 0 \}$ is the non-empty closed, aperiodic  class of communicative states, while $\XX_0 = \{ i \in \XX: \pi_{0, i} = 0 \}$ is the class {\rm (}possibly empty{\rm )} of transient states, for the Markov chain $X_{0, n}$. Note that $\bar{\pi}_0$ still is the stationary distribution for the Markov chain $X_{0, n}$. \vspace{1mm}

We would like also to refer to   paper \cite{MIT}, where one can find alternative upper bounds for the rate of convergence of stationary distributions for perturbed Markov chains and further related references. \vspace{1mm}

{\bf 5.2 Rate of convergence for stationary distributions of singularly perturbed Markov chains with damping component}.
Let now assume  that condition ${\bf B}_1$ holds.

Let us consider matrices, for $j = 1, \ldots, h$ and $n = 0, 1, \ldots$,

\begin{equation}\label{matrad}
\mathbf{P}^n_{0, j} = \| p^{(j)}_{0, rk}(n) \|_{r, k \in \XX^{(j)}}.
\end{equation}

Note that, for $j = 1, \ldots, h$,  probabilities $p^{(j)}_{0, rk}(n) = p_{0, rk}(n), r, k \in \XX^{(j)}, n \geq 0$, since $\XX^{(j)}, j = 1, \ldots, h$ are closed classes of states.

The reduced Markov chain $X^{(j)}_{0, n}$ with the phase space $ \XX^{(j)}$ and the matrix of transition probabilities $\mathbf{P}_{0, j}$ is, for every $j = 1, \ldots, h$, exponentially ergodic and the following estimates take place, for $k \in \XX^{(j)}, j = 1, \ldots, h$ and $n = 0, 1, \ldots$,
\begin{equation}\label{ester}
\max_{r, k \in \XX^{(j)}} |p^{(j)}_{0, rk}(n) - \pi^{(j)}_{0, k} | \leq C_{j} \lambda^n_{j},
\end{equation}
with some constants $C_j = C_j(\mathbf{P}_{0, j}) \in [0, \infty), \lambda_j = \lambda_j(\mathbf{P}_{0, j}) \in [0, 1), j = 1, \ldots, h$ and distributions
$\bar{\pi}^{(j)}_0 = \langle \pi^{(j)}_{0, k}, k\in \XX^{(j)} \rangle, j = 1, \ldots, h$, with all positive components.

Obviously, inequalities (\ref{ester}) imply that $p^{(j)}_{0, rk}(n)  \to  \pi^{(j)}_{0, k}$ as $n \to \infty$, for $r, k \in \XX^{(j)}, j = 1, \ldots, h$. Thus, distribution $\bar{\pi}^{(j)}_0$ is the stationary distribution for the Markov chain $X^{(j)}_{0, n}$, for every $j = 1, \ldots, h$.

As has been mentioned above the role of $\lambda_j$ can play, for every $j = 1, \ldots, h$,  the absolute value of  the second (by absolute value),  eigenvalue for matrix $\mathbf{P}_{0, j}$, and  $C_{j}$ is the constant,  which can be computed using the algorithm described in the above mentioned book \cite{Fel}.

Let us denote,
\begin{equation}\label{erten}
\lambda = \max_{1 \leq j  \leq h} \lambda_j, \ C = \max_{1 \leq j  \leq h} C_j.
\end{equation}

Condition ${\bf B}_1$  is, in fact,  equivalent to the following condition:
\begin{itemize}
\item [${\bf B}_2$:] The phase space $\XX = \cup_{j = 1}^h \XX^{(j)}$, where: (a) $X^{(j)}, j = 0, \ldots, h$ are non-intersecting subsets of $\XX$, (b)
$X^{(j)}, j = 1, \ldots, h$ are non-empty, closed classes of  states for the Markov chain $X_{0, n}$ such that inequalities (\ref{ester}) hold.
\end{itemize}

Indeed, condition ${\bf B}_2$ implies that probabilities $p^{(j)}_{0, rk}(n) > 0, r, k \in \XX^{(j)}, j = 1, \ldots, h$ for all large enough $n$. This implies that $\XX^{(j)}, j = 1, \ldots, h$ are closed, aperiodic classes of communicative states.  Also, inequalities (\ref{ester}) imply that $p^{(j)}_{0, rk}(n)  \to  \pi^{(j)}_{0, k}$ as $n \to \infty$, for $r, k \in \XX^{(j)}, j = 1, \ldots, h$, and, thus, $\bar{\pi}^{(j)}_0 = \langle \pi^{(j)}_{0, k}, k \in \XX^{(j)} \rangle$ is the stationary distribution for the Markov chain $X^{(j)}_{0, n}$, for every $j = 1, \ldots, h$.  \vspace{1mm}

{\bf Theorem 2}. {\em Let condition ${\bf B}_2$ holds. Then the following relation holds, for $k \in \XX$,}
\begin{equation}\label{estima}
| \pi_{\e, k}  - \pi_{0, \bar{d}, k} |  \leq \e ( | d_k - \pi_{0, \bar{d}, k}|  + \frac{C \lambda}{1 - \lambda}).
\end{equation}

The proof of Theorem 2 is analogous to the proof of Theorem 1. \\

{\bf 6. Asymptotic expansions for perturbed stationary \\ \makebox[10mm]{} distributions} \\

In this section, we present asymptotic expansions for stationary distributions of Markov chains with damping component with respect to damping parameter. \vspace{1mm}

{\bf 6.1 Asymptotic expansions for stationary distributions of regularly perturbed Markov chains with damping component}.
Let us get some asymptotic expansions for perturbed stationary distributions in the case, where  condition ${\bf A}_1$ holds.

According to the Perron-Frobenious theorem, in this case, the  eigenvalues $\rho_1, \ldots$, $\rho_m$ of the stochastic matrix $\mathbf{P}_0$ satisfy
the following condition:
\begin{itemize}
\item [${\bf A}_3$:] $\rho_1 = 1 >   |\rho_2| \geq \cdots  \geq |\rho_m|$.
\end{itemize}

Note that some of eigenvalues $\rho_2, \ldots, \rho_m$ can be complex numbers.

As well known, condition ${\bf A}_3$ implies that the following eigenvalues decomposition representation takes place, for $i, j \in \XX$ and $n \geq 1$,
\begin{equation}\label{decom}
p_{0, ij}(n) = \pi_{0, j} + \rho^n_2 \pi_{0, ij}[2] + \cdots + \rho^n_m \pi_{0, ij}[m],
\end{equation}
where: (a) $\bar{\pi}_0 = \langle \pi_{0, j}, j \in \XX \rangle$ is distribution with all positive component, (b) $\pi_{0, ij}[l], i, j \in \XX, l = 2, \ldots, m$
are some real or complex valued coefficients, (c) $\pi_{0, ij}[k] =  \pi_{0, ij}[l], i, j \in \XX$,  if
$\rho_k = \rho_l$, for some $2 \leq k, l \leq m$.

In is useful to note that the eigenvalues decomposition representation can be  rewritten in the alternative form based on
distinct eigenvalues for matrix $\mathbf{P}_0$.
Let $\bar{\rho}_k, k = 1, \ldots, \bar{m}$ be the above mentioned distinct eigenvalues, $\XX_k = \{1 \leq l \leq m: \rho_l = \bar{\rho}_k \}$, and $v_k$ be the index of multiplicity of the eigenvalue $\bar{\rho}_k$ (the number of items in  set $\XX_k$), for $k = 1, \ldots, \bar{m}$. Also, let us assume that $\bar{\rho}_1 = 1  >   |\bar{\rho}_2| \geq \cdots  \geq |\bar{\rho}_{\bar{m}}|$,
and denote $\bar{\pi}_{0, ij}[k] = v_k \pi_{0, ij}[l_k], i, j \in \XX, l_k \in \XX_k, k =2, \ldots, \bar{m}$ (these coefficients
do not depend on the choice of  $l_k \in \XX_k, k = 1, \ldots, \bar{m}$).

Then, the eigenvalues decomposition representation (\ref{decom}) can be rewritten in the following form, for $i, j \in \XX$ and $n \geq 1$,
\begin{equation}\label{decomba}
p_{0, ij}(n) = \pi_{0, j} + \bar{\rho}^n_2 \bar{\pi}_{0, ij}[2] + \cdots +  \bar{\rho}^n_{\bar{m}} \bar{\pi}_{0, ij}[\bar{m}].
\end{equation}

Obviously, relation (\ref{decom}) implies that probabilities $p_{ij}(n) \to \pi_{0, j}$ as $n \to \infty$, for $i, j \in \XX$. Thus,  $\bar{\pi}_0$ is the stationary distribution for the Markov chain $X_{0, n}$.

In fact, condition ${\bf A}_3$ is equivalent to condition ${\bf A}_1$.

Indeed, the convergence  relation, $p_{0, ij}(n) \to \pi_{0, j}$ as $n \to \infty$, for $i, j \in \XX$, implies that   $p_{0, ij}(n) > 0, i, j \in \XX$, for  all large
enough $n$. This implies that $\XX$ is one aperiodic class of communicative states.

We refer to book \cite{Fel}, where one can find the description of effective algorithm for finding matrices $\mathbf{\Pi}_l = \| \pi_{0, ij}[l]  \|, l = 2, \ldots, m$.

Relation (\ref{decom}) implies the following relation holds, for  $j \in \XX$ and $n \geq 1$,
\begin{equation}\label{decoma}
p_{0, \bar{d}, j}(n) = \pi_{0, j} + \rho^n_2 \pi_{0, \bar{d}, j}[2] + \cdots + \rho^n_m \pi_{0,\bar{d}, j}[m],
\end{equation}
where, for $j \in \XX, l = 2, \ldots, m$,
\begin{equation}\label{goptrewa}
\pi_{0, \bar{d}, j}[l]  = \sum_{i \in \XX} d_i \pi_{0, i j}[l].
\end{equation}

Let also define coefficients, for  $j \in \XX, n \geq 1$,
\begin{equation}\label{last}
\tilde{\pi}_{0, \bar{d}, j}[n] = \left\{
\begin{array}{lll}
d_j - \pi_{0, j} + \sum_{l = 2}^m  \pi_{0, \bar{d}, j}[l]  \frac{\rho_l}{1 - \rho_l} & \text{for} \ n = 1, \vspace{2mm}  \\
(-1)^{n -1} \sum_{l = 2}^m  \pi_{0, \bar{d}, j}[l]   \frac{\rho_l^{n-1}}{(1 - \rho_l)^{n }} & \text{for} \ n  > 1.
\end{array}
\right.
\end{equation}

Below,  symbol $O(\e^n)$ is  used for quantitiies such that  $O(\e^n)/\e^n$ is
bounded  as function of $\e \in (0, 1]$.

The following theorem takes place. \vspace{1mm}

{\bf Theorem 3}. {\em Let condition ${\bf A}_3$ holds. Then,  the following asymptotic expansions take place for every $j \in \XX$ and $n \geq 1$,}
\begin{equation}\label{srewavabe}
\pi_{\e, j} = \pi_{0, j} +  \tilde{\pi}_{0, \bar{d}, j}[1] \e + \cdots +  \tilde{\pi}_{0, \bar{d},  j}[n] \e^{n}  + O(\e^{n+1}).
\end{equation}

{\bf Proof}. Relations (\ref{gotyrtbaramoto}) and  (\ref{decoma}) imply that the following relation holds, for  $j \in \XX$,
\begin{align}\label{srewa}
\pi_{\e, j} & = \e \sum_{n = 0}^{\infty}  p_{0, \bar{d}, j}(n)(1- \e)^n \vspace{3mm} \nonumber \\
& = \e d_j + \e \sum_{n = 1}^{\infty} (\pi_{0, j} +  \sum_{l = 2}^m \rho^n_l \pi_{0, \bar{d}, j}[l] )(1- \e)^n \vspace{3mm} \nonumber \\
& = \pi_{0, j}  + \e(d_j - \pi_{0, j}) + \sum_{l = 2}^m \pi_{0, \bar{d}, j}[l]  \e \sum_{n = 1}^{\infty} \rho^n_l  (1- \e)^n \vspace{3mm} \nonumber \\
& = \pi_{0, j} + \e(d_j - \pi_{0, j}) + \sum_{l = 2}^m  \pi_{0, \bar{d}, j}[l]   \frac{\rho_l \e (1- \e)}{1 - \rho_l(1-\e)} \vspace{3mm} \nonumber \\
& = \pi_{0, j} + \e(d_j - \pi_{0, j})   + \sum_{l = 2}^m  \pi_{0, \bar{d}, j}[l]  \rho_l  \e (1 - \e) (1 - \rho_l + \rho_l  \e)^{-1}.
\end{align}

Functions $(a + b\e)^{-1}, \e \in [0, 1]$ and $b \e(1 - \e)(a + b\e)^{-1}, \e \in [0, 1]$ admit, for any complex numbers $a \neq 0$ and $b$,  the following Taylor asymptotic
expansions, for every $n \geq 1$ and $\e \to 0$,
\begin{align}\label{expa}
(a + b\e)^{-1} & =  a^{-1} -  a^{-2}b \e +  a^{-3}b^2\e^2 \nonumber \\
& \quad + \cdots + (-1)^{n} a^{- (n + 1)}b^n \e^n + O(\e^{n+1}),
\end{align}
and
\begin{align}\label{expana}
& b\e (1 - \e) (a + b\e)^{-1}  =  a^{-1}b\e -  a^{-1}b (1 + a^{-1}b) \e^2 +  a^{-2}b^2(1 + a^{-1}b) \e^3 \nonumber \\
& \quad \quad \quad + \cdots + (-1)^{n-1} a^{- (n -1)}b^{n-1}(1 + a^{-1}b) \e^{n} + O(\e^{n+1}).
\end{align}

Relations  (\ref{srewa}) -- (\ref{expana}) let us write down the following Taylor asymptotic
expansions for stationary probabilities $\pi_{\e, j}, j \in \XX$, for every
$n \geq 1$ and $\e \to 0$,
\begin{align}\label{srewava}
\pi_{\e, j} & = \pi_{0, j} + \e(d_j - \pi_{0, j})  + \sum_{l = 2}^m  \pi_{0, \bar{d}, j}[l] \rho_l \e(1 - \e) (1 - \rho_l + \rho_l \e)^{-1} \vspace{3mm} \nonumber \\
& = \pi_{0, j} +  \tilde{\pi}_{0, \bar{d}, j}[1] \e + \cdots +  \tilde{\pi}_{0, \bar{d}, j}[n] \e^{n}  + O(\e^{n+1}).
\end{align}

The proof is complete. $\Box$ \vspace{1mm}

It worth noting that some of eigenvalues $\rho_l$ and coefficients $\pi_{0, ij}[r]$ can be complex numbers. Despite of this, coefficients 
$\tilde{\pi}_{0, \bar{d},  j}[n], n \geq 1$ in the expansions given in relation {\rm (\ref{srewavabe})} are real numbers.

Indeed, $\pi_{\e, j}$  is a positive number,  for every
$\e \in [0, 1]$. Relation {\rm (\ref{srewava})} implies
that  $(\pi_{\e, j}  - \pi_{0, j}) \e^{-1} \to \tilde{\pi}_{0, \bar{d}, j}[1]$ as $\e \to 0$. Thus, $\tilde{\pi}_{0, \bar{d}, j}[1]$ is a real number. In this way, the above proposition can be proved for all coefficients in expansions {\rm (\ref{srewavabe})}. This implies that the remainders of these expansions  $O(\e^{n+1})$ also are real-valued functions of $\e$. 

Moreover, since $\bar{\pi}_\e = \langle \pi_{\e, j}, j \in \XX \rangle, \e \in (0, 1]$ and  $\bar{\pi}_{0} = \langle \pi_{0, j}, j \in \XX \rangle$ are probability distributions, the following equalities connect coefficients in the asymptotic expansions (\ref{srewavabe}), 
$\sum_{j \in \XX} \tilde{\pi}_{0, \bar{d}, j}[n] = 0$, for $n \geq 1$.

\vspace{1mm}

{\bf 6.2 Asymptotic expansions for stationary distributions of singularly perturbed Markov chains with damping component}.
 Let now consider the case, where   condition ${\bf B}_1$ holds.
We can assume that class $\XX^{(j)}$ includes $m_j$ states, for $j = 1, \ldots, h$, where  $m_j > 0, j = 1, \ldots, h$.

Let us denote by $\rho_{j, 1}, \ldots$, $\rho_{j, m_j}$ the eigenvalues  of the stochastic matrices
$\mathbf{P}_{0, j}$, $j = 1, \ldots, h$. We can assume that these eigenvalues  are ordered by absolute values , i.e. $|\rho_{j, 1}| \geq    |\rho_{j, 2}| \geq \cdots  \geq |\rho_{j, m_j}|$, for $j = 1, \ldots, h$.

Condition ${\bf B}_1$ is, in fact, equivalent to the following condition:
\begin{itemize}
\item [${\bf B}_3$:] The phase space $\XX = \cup_{j = 1}^h \XX^{(j)}$, where: (a) $X^{(j)}, j = 1, \ldots, h$ are non-intersecting subsets of $\XX$, (b)
$X^{(j)}, j = 1, \ldots, h$ are non-empty, closed classes of  states for the Markov chain $X_{0, n}$ such that  inequalities   $\rho_{j, 1} = 1 >   |\rho_{j, 2}| \geq \cdots  \geq |\rho_{j, m_j}|, \,  j = 1, \ldots, h$,   hold.
\end{itemize}

Condition ${\bf B}_1$ implies that condition ${\bf B}_3$ holds.  This follows from the Perron-Frobenius theorem.

 Condition ${\bf B}_3$ imply that the following eigenvalues decomposition representations take place, for $r, k \in \XX^{(i)}, j = 1, \ldots, h$ and $n \geq 1$,
\begin{equation}\label{decombanas}
p^{(j)}_{0, rk}(n) = \pi^{(j)}_{0, k} + \rho_{j, 2}^n \pi^{(j)}_{0, rk}[2] + \cdots + \rho_{j, m_j}^n \pi^{(j)}_{0, rk}[m_j],
\end{equation}
where:  (a) $\bar{\pi}^{(j)}_0 = \langle \pi^{(j)}_{0, k}, k \in \XX^{(j)} \rangle$ is a distribution with all positive component, for $j = 1, \ldots, h$,  (b) $\pi^{(j)}_{0, rk}[l], r, k \in \XX^{(j)}, l= 2, \ldots, m_j, j = 1, \ldots, h$  are some real or complex valued coefficients  satisfying equalities $\pi^{(j)} _{0, rk}[s] =  \pi^{(j)} _{0, rk}[t], r, k \in \XX^{(j)}$,  if
$\rho_{j, s} = \rho_{j, t}$, for some $2 \leq s, t \leq m_j, j = 1, \ldots, h$.

Obviously, relation (\ref{decomba}) implies that probabilities $p^{(j)}_{rk}(n) \to \pi_{0, j}$ as $n \to \infty$, for $r, k \in \XX^{(j)}, j = 1, \ldots, h$. Thus,  $\bar{\pi}^{(j)}_0$ is the stationary distribution of the Markov chain $X^{(j)}_{0, n}$, for $j = 1, \ldots, h$.

In fact, condition ${\bf B}_3$ is equivalent to condition ${\bf B}_1$.

Indeed, as was mentioned above, relation (\ref{decomba}) implies that  probabilities $p^{(j)}_{0, rk}(n) \to \pi^{(j)}_{0, k}$ as $n \to \infty$, for $r, k \in \XX^{(j)}, j = 1, \ldots, h$. Thus,  probabilities    $p^{(j)}_{0, rk}(n) > 0, r, k \in \XX, j = 1, \ldots, h$, for  all large enough $n$. This implies that $\XX^{(j)}$ is, for every $j = 1, \ldots, h$ a closed aperiodic class of communicative states for the Markov chain $X_{0, n}$.
Thus, condition ${\bf B}_3$ implies that condition ${\bf B}_1$ holds.

We  again refer to book \cite{Fel}, where one can find the description of effective algorithm for finding matrices
$\mathbf{\Pi}^{(j)}_l = \| \pi^{(j)}_{0, rk}[l]  \|, l = 2, \ldots, m_j, j = 1, \ldots, h$.

Relation (\ref{decomba}) implies, in this case, the following relation holds, for any  $k \in \XX^{(j)}, j = 1, \ldots, h$ and $n \geq 1$,
\begin{equation}\label{decomane}
p_{0, \bar{d}, k}(n) = \pi_{0, \bar{d}, k} + \rho^n_{j, 2} \pi^{(j)}_{0, \bar{d}, k}[2] + \cdots + \rho^n_{j, m_j} \pi^{(j)}_{0,\bar{d}, k}[m_j],
\end{equation}
where, for $k \in \XX^{(j)},  l = 2, \ldots, m_j,  j = 1, \ldots, h$,
\begin{equation}\label{goptrewane}
\pi^{(j)}_{0, \bar{d}, k}[l]  = \sum_{r \in \XX^{(j)}} d_r \pi^{(j)}_{0, r k}[l].
\end{equation}

Let also define coefficients, for  $k \in \XX^{(j)}, j = 1, \ldots, h, n \geq 1$,
\begin{equation}\label{lastnem}
\tilde{\pi}^{(j)}_{0, \bar{d}, k}[n]  = \left\{
\begin{array}{lll}
d_k - \pi^{(j)}_{0, k} + \sum_{l = 2}^{m_j}  \pi^{(j)}_{0, \bar{d}, k}[l]  \frac{\rho_{j, l}}{1 - \rho_{j, l}} & \text{for} \ n = 1, \vspace{2mm}  \\
(-1)^{n -1}  \sum_{l = 2}^{m_j}  \pi^{(j)}_{0, \bar{d}, k}[l]  \frac{\rho_{j, l}^{n-1}}{(1 - \rho_{j, l})^{n }} & \text{for} \ n  > 1.
\end{array}
\right.
\end{equation}

The following theorem takes place. \vspace{1mm}

{\bf Theorem 4}. {\em Let condition ${\bf B}_3$  holds. Then,  the following asymptotic
expansion take place for every for  $k \in \XX^{(j)}, j = 1, \ldots, h, n \geq 1$,}
\begin{equation}\label{srewavabena}
\pi_{\e, k} = \pi_{0, \bar{d},  k} +  \tilde{\pi}^{(j)}_{0, \bar{d}, k}[1] \e + \cdots +  \tilde{\pi}^{(j)}_{0, \bar{d},  k}[n] \e^{n}  + O(\e^{n+1}).
\end{equation}

The proof is similar with  the proof of Theorem 3.  \\

 {\bf 7. Coupling  for regularly perturbed Markov chains \\ \makebox[10mm]{}  with damping component} \\

In this section, we present  coupling algorithms and get the effective upper bound for the rate of convergence in ergodic theorems for regularly perturbed Markov chains with damping component. \vspace{1mm}

{\bf 7.1 Maximal coupling for discrete distributions}. Let $\bar{p}' = \langle p'_i, i \in \XX \rangle$ and $\bar{p}'' = \langle p''_1, i \in \XX  \rangle$ be two
discrete probability distributions. Let us denote by ${\cal L}[\bar{p}', \bar{p}'']$ the class of   two-dimensional distribution $\bar{P} =   \langle  P_{ij}, (i,j) \in \XX \times \XX \rangle$ which satisfy the following conditions  (a) $P'_i = \sum_{j \in \XX} P_{ij} = p'_i, i \in \XX$; (b) $P''_j = \sum_{i \in \XX} P_{ij} = p''_j, j \in \XX$.

Let us also denote,
\begin{equation}\label{qewr}
Q_{\bar{P}} = \sum_{i \in \XX} P_{ii}
 \end{equation}
 and
 \begin{equation}\label{qewra}
 Q(\bar{p}', \bar{p}'') = \sup_{\bar{P} \in {\cal L}[\bar{p}', \bar{p}'']} Q_{\bar{P}}.
 \end{equation}

 The following lemma presents the well known ``coupling'' result, which variants can be found in \cite{Gri, Kal, Lin, Pit} and \cite{Si4, Si5, Si6, Si7, SPe}. \vspace{1mm}

  {\bf Lemma 6}. {\em There exists the two-dimensional distribution $\bar{P}^*  = \langle P_{ij}^*, i, j \in \XX \rangle  \in  {\cal L}[\bar{p}', \bar{p}'']$ such that:
 \begin{equation}\label{sol}
Q_{\bar{P}^*} =  Q^* = \sum_{i \in \XX} \min(p'_i, p''_i) = Q(\bar{p}', \bar{p}'').
 \end{equation}

The  distribution $\bar{P}^*$ is given by the following relations{\rm :} \vspace{1mm}

 {\bf (i)} If  $Q^* \in (0, 1)$, then
\begin{align}\label{solv}
P^*_{ij} = &  \min(p'_i, p''_j) {\rm I} (i = j) \vspace{1mm} \nonumber \\
& \, + \frac{1}{1 - Q^*} ( p'_i - \min(p'_i, p''_i))( p''_j - \min(p'_j, p''_j)), i, j \in \XX.
 \end{align}

 {\bf (ii)}  $Q^* = 1$ if and only if $ p'_k = p''_k, k \in \XX$ and,  in this case,
 \begin{equation}\label{solva}
 P^*_{ij} = \min(p'_i, p_j'') {\rm I} (i = j), \,  i, j \in \XX. \makebox[45mm]{}
 \end{equation}

  {\bf (iii)} $Q^* = 0$ if and only if $p'_k p''_k = 0, k \in \XX$ and, in this case,
 \begin{equation}\label{solvab}
  P^*_{ij} = p'_i p''_j, i, j \in \XX. \makebox[72mm]{}
 \end{equation}}

 {\bf Proof}. It can be  found in the above mentioned works. In order, to improve self-readability of the present paper, we just give a short sketch of the proof. Obviously, probability $P_{ii} \leq p'_i \wedge p''_i, i \in \XX$, for any two-dimensional distribution $\bar{P} =
\langle  P_{ij}, (i,j) \in \XX \times \XX \rangle \in {\cal L}[\bar{p}', \bar{p}'']$. This relation implies that  $Q_{\bar{P}} \leq Q^* =  \sum_{i \in \XX} p'_i \wedge p''_i$.  This is easily to check that every relation (\ref{solv}), (\ref{solva}), or (\ref{solvab}) defines a two-dimensional distribution $\bar{P}^*$ from the class ${\cal L}[\bar{p}', \bar{p}'']$.  Moreover, the corresponding  quantity $Q_{\bar{P}^*} =  Q^*$. This is obvious for two cases presented in propositions  {\bf (ii)}  and  {\bf (iii)}.    In the first case presented in
proposition  {\bf (i)}, this follows from the relation,  $( p'_i - \min(p'_i, p''_i))( p''_i - \min(p'_i, p''_i)) = 0, i \in \XX$. $\Box$ \vspace{1mm}

 {\bf 7.2  A one-step coupling for regularly perturbed  Markov chains with damping component}.
 Let us recall the matrix of transition probabilities $\mathbf{P}_\e  = \| p_{\e, ij} \|$ and define the quantities,
 \begin{equation}\label{utre}
Q_{\e, ij} = \sum_{r \in \XX}\min(p_{\e, ir}, p_{\e, jr}), i, j \in \XX.
\end{equation}

Let now use the ``one-step''  coupling algorithm for construction a ``coupling''  Markov chain $Z_{\e, n} =  (X'_{\e, n}, X''_{\e, n}), n = 0, 1, \ldots$, with:

{\bf (i)} the phase space $\ZZ = \XX \times \XX$;

{\bf (ii)} the initial distribution $\bar{P}_{\e} = \langle P_{\e, ij} , (i, j) \in \ZZ \rangle$ constructing according to relation (\ref{solv}), (\ref{solva}), or
(\ref{solvab}) for distributions $\bar{p}' = \bar{p} = \langle p_i, i \in \XX \rangle$ and
$\bar{p}'' = \bar{\pi}_{\e} = \langle \pi_{\e, i}, i \in \XX \rangle$;

{\bf (iii)} the transition probabilities $P_{\e, ij, rk}$ defined by the following relations, for $(i, j), (r, k) \in \ZZ$:
 \vspace{1mm}

\noindent {\bf (a)} If $Q_{\e, ij} \in (0, 1)$, then,
 \begin{align}\label{gotyrbedga}
 P_{\e, ij, rk} &  =  \PP \{ X'_{\e, 1} = k, X''_{\e, 1} = r / X'_{\e, 0} = i, X''_{\e, 0} = j \} \vspace{2mm} \nonumber \\
& = \min(p_{\e, ir}, p_{\e, jk}) {\rm I} (r = k) \vspace{2mm} \nonumber \\
& \quad  + \frac{1}{1 - Q_{\e, ij}} ( p_{\e, ir} - \min(p_{\e, ir}, p_{\e, jr}))( p_{\e, jk} - \min(p_{\e, ik}, p_{\e, jk})).
 \end{align}

\noindent {\bf (b)}  If $Q_{\e, ij} = 1$, then $p_{\e, ir} = p_{\e, jr}, r \in \XX$ and,
  \begin{equation}\label{gotyrbedgab}
 P_{\e, ij, rk} = \min(p_{\e, ir}, p_{\e, jk}) {\rm I} (r = k). \makebox[62mm]{}
 \end{equation}

\noindent  {\bf (c)}  If $Q_{\e, ij} = 0$, then $p_{\e, ir} p_{\e, jr} = 0, r \in \XX$ and,
  \begin{equation}\label{gotyrbedgabe}
  P_{\e, ij, rk}  = p_{\e, ir} p_{\e, jk}. \makebox[88mm]{}
 \end{equation}

The above construction of coupling Markov chain and the following lemma originate from works \cite{Gri} and \cite{Pit} and  plays an important role in what follows.

\vspace{1mm}

{\bf Lemma 7}. {\em Let $Z_{\e, n} =  (X'_{\e, n}, X''_{\e, n}), n =0, 1, \ldots$ be a homogeneous Markov chain with the phase space $\ZZ = \XX \times \XX$, the initial distribution $\bar{P_\e}$ and transition probabilities given by relations {\rm (\ref{gotyrbedga}) -- (\ref{gotyrbedgabe})}. Then{\rm :} \vspace{1mm}

 {\bf (i)} The first component, $X'_{\e, n}, n = 0, 1, \ldots$, is  a homogeneous Markov chain with the phase space $\XX$, the initial distribution $\bar{p}$ and the matrix of transition probabilities $\mathbf{P}_\e$.

{\bf (ii)} The second component  $X''_{\e, n}, n = 0, 1, \ldots$  is a homogeneous Markov chain with the phase space $\XX$, the initial distribution $\bar{\pi}_\e$ and the matrix of transition probabilities $\mathbf{P}_\e$.

{\bf (iii)} The set $\ZZ_0 = \{ (i, i), i \in \XX \}$ is an absorbing set for the  Markov chain $Z_{\e, n}$, i.e.,  probabilities $P_{\e, ii, rk} = 0$, for $i, r, k  \in \XX, r \neq k$.} \vspace{1mm}

{\bf Proof}. Variants of the proof can be  found in the above mentioned works. In order, to improve self-readability of the present paper, we just give a short sketch of the proof. Let us consider the case, where transition probabilities $P_{\e, ij, rk}$ are given by relation (\ref{gotyrbedga}). According to Lemma 6, the following relation takes place, for $i, j, r \in \XX$,
\begin{align}\label{hutras}
& \sum_{k \in \XX}  P_{\e, ij, rk}  =  \sum_{k \in \XX} \min(p_{\e, ir}, p_{\e, jk}){\rm I} (r = k)  \vspace{2mm} \nonumber \\
& \quad  \quad + \frac{1}{1 - Q_{\e, ij}} ( p_{\e, ir} - \min(p_{\e, ir}, p_{\e, jr}))
\sum_{k \in \XX}( p_{\e, jk} - \min(p_{\e, ik}, p_{\e, jk})) \makebox[5mm]{} \vspace{2mm} \nonumber \\
& \quad = \min(p_{\e, ir}, p_{\e, jr}) + p_{\e, ir} - \min(p_{\e, ir}, p_{\e, jr}) = p_{\e, ir}.
\end{align}

The proof of relation (\ref{hutras}) in the cases  where transition probabilities $P_{\e, ij, rk}$ are given by relation (\ref{gotyrbedgab}) or
(\ref{gotyrbedgabe}) is trivial.

Using the assumption that $Z_{\e, n}$ is a Markov chain and  relation (\ref{hutras}), we get the following relation, for any chain of states $i_0, \ldots, i_n \in \XX, n \geq 0$,
\begin{align}\label{hutrasa}
& \PP \{ X'_{\e, l} = i_l, l = 0, \ldots, n \} = \sum_{j_0, \ldots, j_n \in \XX} \PP \{ X'_{\e, l} = i_l, X''_{\e, l} = j_l, l = 0, \ldots, n \}  \makebox[5mm]{} \vspace{2mm} \nonumber \\
& \quad  \quad = \sum_{j_0 \in \XX} p_{i_0} \pi_{\e, j_0}  \sum_{j_1 \in \XX} P_{\e, i_0 j_0, i_1 j_1} \cdots
\sum_{j_{n} \in \XX} P_{\e, i_{n-1} j_{n-1}, i_n j_n} \vspace{2mm} \nonumber \\
& \quad  \quad =  \sum_{j_0 \in \XX} p_{i_0} \pi_{\e, j_0}  \sum_{j_1 \in \XX} P_{\e, i_0 j_0, i_1 j_1} \cdots
\sum_{j_{n-1} \in \XX} P_{\e, i_{n-2} j_{n-2}, i_{n-1} j_{n-1}} \cdot p_{\e, i_{n-1}, i_n} \vspace{2mm} \nonumber \\
& \quad \quad  = \cdots = \sum_{j_0 \in \XX} p_{i_0} \pi_{\e, j_0} p_{\e, i_{0}, i_1} \cdots p_{\e, i_{n-1}, i_n} \vspace{2mm} \nonumber \\
& \quad  \quad  = p_{i_0} p_{\e, i_{0}, i_1} \cdots p_{\e, i_{n-1}, i_n}.
\end{align}

Relation (\ref{hutrasa}) proves proposition {\bf (i)}. The proof of proposition {\bf (ii)} is analogous.

According to relation (\ref{gotyrbedgab}), the quantities $Q_{\e, ii} = 1, i \in \XX$. This proves proposition {\bf (iii)}. $\Box$ \vspace{1mm}

Let $\mathbf{A} = \| a_{ij} \|$ be a $m\times m$ a matrix with real-valued elements. Let us introduce functional,
\begin{equation}\label{hopiun}
Q(\mathbf{A}) = \min_{1 \leq i, j \leq m} \sum_{k = 1}^m a_{ik} \wedge a_{jk}.
\end{equation}

The following simple lemma presents some basic properties of functional $Q(\mathbf{A})$. \vspace{1mm}

{\bf Lemma 8}. {\em  Functional $Q(\mathbf{A})$ possesses the following properties{\rm :} {\rm (a)} $Q(a\mathbf{A}) = a Q(\mathbf{A})$, for any $a \geq 0${\rm ;}  {\rm (b)} $Q(\mathbf{A}) \geq  a_1Q(\mathbf{A}_1) + \cdots + a_n Q(\mathbf{A}_n)$, for  any $m\times m$  matrices
$\mathbf{A}_1, \ldots, \mathbf{A}_n$ with real-valued elements, numbers $a_1, \ldots, a_n$ $\geq 0, a_1 \cdots + a_n = 1$, and matrix $\mathbf{A} = a_1\mathbf{A}_1 + \cdots +  a_n \mathbf{A}_n$,  for $n \geq 2${\rm ;}     {\rm (c)}  $Q(\mathbf{A}) \in [0, 1]$,
for any stochastic matrix $\mathbf{A}${\rm ;}   {\rm (d)} $Q(\mathbf{A}) = 1$,  for any $m\times m$  stochastic damping type matrix $\mathbf{A} = \| a_{ij} \|$,  with elements $a_{ij} = a_j \geq 0, i, j = 1, \ldots, m$.} \vspace{1mm}

Lemma 8 implies that the following inequality holds,
\begin{align}\label{bedga}
 Q(\mathbf{P}_\e) & \geq  (1 - \e) Q(\mathbf{P}_0) + \e Q(\mathbf{D}) \vspace{1mm} \nonumber \\
 & = (1- \e) Q(\mathbf{P}_0) + \e,
 \end{align}

Let us also introduce quantities,
\begin{equation}\label{inita}
Q_{\e, \bar{p}} =  Q(\bar{p}, \bar{\pi}_\e) =   \sum_{i \in \XX} \min(p_i, \pi_{\e, i}).
\end{equation}

The given below Theorems 5 and 6 present effective upper bounds for the rate of convergence  in the individual ergodic theorem for perturbed  Markov chain with damping component based on corresponding general coupling results for Markov chains given in
\cite{Gri, Lin, Pit}. Theorems 5 and 6  specify and  detail coupling upper bounds for the rate of convergence
for Markov chains with damping component.

Note that condition ${\bf A}_1$ is not required in Theorem 5 formulated below.

Also, we count $(1 - Q(\mathbf{P}_0))^0 = 1$, if $1 - Q(\mathbf{P}_0) = 0$.
\vspace{1mm}

{\bf Theorem 5}. {\em The following relation takes place, for every $\e \in (0, 1]$ and $\bar{p} \in {\cal P}_m, j \in \XX, n \geq 0$,}
\begin{equation}\label{hoptre}
| p_{\e, \bar{p}, j}(n) - \pi_{\e, j} | \leq ( 1 - Q_{\e, \bar{p}}) (1 - Q(\mathbf{P}_0))^n (1 - \e)^n,
\end{equation}

{\bf Proof}. Obviously, for $\bar{p} \in {\cal P}_m, j \in \XX, n \geq 0$,
\begin{equation}\label{stata}
\PP\{ X'_{\e, n} = j \} = p_{\e, \bar{p}, j}(n).
\end{equation}

Since, the initial distribution of Markov chain $X''_{\e, n}$ coincides with its stationary distribution, this Markov chain is a stationary random sequence and, thus, for $j \in \XX, n \geq 0$,
\begin{equation}\label{statop}
\PP \{ X''_{\e, n} = j \} = \pi_{\e, j}.
\end{equation}

Let now define the hitting (coupling) time,
\begin{equation}\label{stat}
T_\e = \min(n \geq 0: X'_{\e, n} = X''_{\e, n} ) = \min (n \geq 0: Z_{\e, n} \in \ZZ_0).
\end{equation}

Since $\ZZ_0$ is an absorbing set for the Markov chain $Z_{\e, n}$,  the following relation holds,
\begin{equation}\label{statas}
\PP \{ Z_{\e, n} \in \ZZ_0,  n \geq T_\e \} = 1.
\end{equation}

Using the above remarks, we get the following relation, for $j \in \XX, n \geq 0$,
\begin{align}\label{relat}
| p_{\e, \bar{p}, j}(n) - \pi_{\e, j} |  &  = | \PP \{ X'_{\e, n} = j \} -  \PP \{ X''_{\e, n} = j \} |  \vspace{1mm} \nonumber \\
& =   | \PP \{ X'_{\e, n} = j,  X''_{\e, n} \neq j \} -  \PP \{ X'_{\e, n} \neq  j,  X''_{\e, n} = j \} | \vspace{1mm} \nonumber \\
& \leq  \PP \{ X'_{\e, n} = j,  X''_{\e, n} \neq j \} + \PP \{ X'_{\e, n} \neq  j,  X''_{\e, n} = j \}  \vspace{1mm} \nonumber \\
& \leq \PP \{ T_\e > n \}.
\end{align}

Using Lemma 6, we get,  $\bar{p} \in {\cal P}_m, j \in \XX$,
\begin{align}\label{serof}
| p_{\e, \bar{p}, j}(0) - \pi_{\e, j} | & \leq  \PP \{ X'_{\e, 0} = j,  X''_{\e, 0} \neq j \} + \PP \{ X'_{\e, 0} \neq  j,  X''_{\e, 0} = j \}
\vspace{1mm} \nonumber \\
& \leq \PP \{ T_\e  >  0 \} =  1 -  Q_{\e, \bar{p}}.
\end{align}

Also, by continuing the inequalities (\ref{relat}) we get, for  $\bar{p} \in {\cal P}_m, j \in \XX, n \geq 1$,
\begin{align}\label{relatas}
| p_{\e, \bar{p}, j}(n) - \pi_{\e, j} | & \leq  \PP \{ T_\e > n \} \vspace{2mm} \nonumber \\
 & =    \sum_{i, j \in \XX} \PP \{  T_\e > n -1,  X'_{\e, n-1} = i, X''_{\e, n-1} = j \}  \vspace{2mm} \nonumber \\
& \quad \quad \quad \times  \PP \{ X'_{\e, n} \neq X'_{\e, n} / X'_{\e, n -1} = i,  X'_{\e, n -1} = j \} \vspace{2mm} \nonumber \\
& = \sum_{i, j \in \XX} \PP \{ T_\e > n -1, X'_{\e, n-1} = i, X''_{\e, n-1} = j \} (1 -  Q_{\e, ij})
\vspace{1mm} \nonumber \\
& \leq \PP \{  T_\e > n -1 \} (1 - Q(\mathbf{P}_\e)) \vspace{2mm} \nonumber \\
& \leq  \cdots \leq \PP \{  T_\e > 0 \} (1 - Q(\mathbf{P}_\e))^n   \vspace{2mm} \nonumber \\
& = (1 - Q_{\e, \bar{p}})(1 - Q(\mathbf{P}_0))^n (1 - \e)^n .
\end{align}

The proof is complete. $\Box$ \vspace{1mm}

{\bf Remark 5}. The upper bounds given in relation {\rm (\ref{hoptre})} become better if quantities $1 - Q_{\e, \bar{p}}$, $1 - Q(\mathbf{P}_0)$ and  $1 - \e$ take smaller values. The factor $1 - Q_{\e, \bar{p}}$, is determined by a prior information about the stationary probabilities. It takes smaller values if one can choose initial distribution $\bar{p}$ with smaller deviation from the stationary distribution $\bar{\pi}_\e$.
Relation {\rm (\ref{hoptre})} gives an effective upper bounds for the rate of convergence in the corresponding individual ergodic theorem even in the case where
factor $1 - Q(\mathbf{P}_0) = 1$. \vspace{1mm}

It worth noting that the upper bound $(1- \e)^n$ on the the right hand side of inequality  (\ref{hoptre}) have been given for Markov
chains with a general phase  and damping component,  in the recent paper  \cite{APZ}.

\vspace{1mm}

 As it was shown in \cite{Pit}, the following formula takes place,
\begin{align}\label{formad}
\Delta_1(\mathbf{P}_0) & = 1 - Q(\mathbf{P}_0) = \max_{i, j \in \XX} (1 - \sum_{k \in \XX} p_{0, ik} \wedge  p_{0, jk}) \vspace{2mm} \nonumber \\
& = \frac{1}{2} \max_{i, j \in \XX} \big( \sum_{k \in \XX} (p_{0, ik} - p_{0, ik} \wedge  p_{0, jk}) + \sum_{k \in \XX} (p_{0, jk} - p_{0, ik} \wedge  p_{0, jk})\big) \vspace{2mm} \nonumber \\
& = \frac{1}{2} \max_{i, j \in \XX} \big( \sum_{k: p_{0, ik} \geq p_{0, jk}} (p_{0, ik} - p_{0, jk}) + \sum_{k: p_{0, jk} \geq p_{0, ik}} (p_{0, jk} - p_{0, ik})\big)
\vspace{2mm} \nonumber \\
& = \frac{1}{2} \max_{i, j \in \XX} \sum_{k \in \XX} | p_{0, ik} -  p_{0, jk}|.
\end{align}

The quantity $\Delta_1(\mathbf{P}_0)$ is known  as an ergodicity coefficient. It  determines ergodic
properties of the Markov chain $X_{0, n}$.

This Markov chain is ergodic with an exponential rate of convergence in the corresponding ergodic theorem, if the following condition
holds:
\begin{itemize}
\item [${\bf C}^{(1)}$:] $\Delta_1(\mathbf{P}_0) \in [0, 1)$.
\end{itemize}

Condition ${\bf C}^{(1)}$ is sufficient for holding the weaken variant of condition ${\bf A}_1$ mentioned in Remark 4.

Indeed, probabilities $p_{\e, \bar{p}, j}(n) \to p_{0, \bar{p}, j}(n)$ as $\e \to 0$, for any $j \in \XX, n \geq 0$. Since stationary probabilities $\pi_{\e, j} \in [0, 1], j \in \XX$, any sequence $0 < \e_n \to 0$ as $n \to \infty$ contain a subsequence $0 < \e_{n_l} \to 0$ as $l \to \infty$ such that $\pi_{\e_{n_l}, j} \to  \pi_{0, j}$ as $l \to \infty$, for $j \in \XX$. By passing $\e \to 0$ in the inequality (\ref{hoptre}), we get the
following relation holding  for  $\bar{p} \in {\cal P}_m, j \in \XX, n \geq 0$,
\begin{equation}\label{hubasqal}
| p_{0, \bar{p}, j}(n) - \pi_{0, j} | \leq
( 1 - Q_{0, \bar{p}}) \Delta_1(\mathbf{P}_0)^n,
\end{equation}
where one should count $\Delta_1(\mathbf{P}_0)^0 = 1$, if $\Delta_1(\mathbf{P}_0) = 0$.

Relation (\ref{hubasqa}) obviously implies that $p_{0, \bar{p}, j}(n) \to \pi_{0, j}$ as $n \to \infty$, for $j \in \XX$. Thus limits $\pi_{0, j}, j \in \XX$ are the same for any subsequences $\e_n$ and $\e_{n_l}$ and, thus,
stationary probabilities $\pi_{\e, j} \to \pi_{0, j}$ as $\e \to 0$, for $j \in \XX$.  By passing $\e \to 0$ in the equations given in relation (\ref{ergom}), we get that limits $\pi_{0, j}, j \in \XX$ satisfy the system of linear equations
(\ref{ergsys}) and, thus, $\bar{\pi}_0 = \langle \pi_{0, j}, j \in \XX \rangle$ is the stationary distribution for the Markov chain $X_{0, n}$. Some components of this stationary distribution can equal $0$. In this case, set $X^{(1)} =
\{ j \in \XX: \pi_{0, j} > 0 \}$ is a closed, aperiodic class of communicative states, while set $X^{(0)} = \{ j \in \XX: \pi_{0, j} = 0 \}$ is the class of transient states, for the Markov chain $X_{0, n}$.

If the stationary distribution  $\bar{\pi}_0 = \langle \pi_{0, j}, j \in \XX \rangle$ is positive, then  $X^{(0)}  = \emptyset$. In this case, condition ${\bf C}^{(1)}$ is sufficient for holding of condition ${\bf A}_1$.

Relation (\ref{hubasqal}) implies that the Markov chain  $X_{0, n}$ is ergodic with an exponential rate of convergence in the corresponding ergodic theorem, if condition ${\bf C}^{(1)}$ holds. \vspace{1mm}

{\bf 7.3  A  multi-step coupling for regularly perturbed  Markov chains with damping component}.  The problem with the ergodicity
coefficient $\Delta_1(\mathbf{P}_0)$ is that it can take value $1$ even in the case where condition  ${\bf A}_1$ holds. This happens if there exists two states $i, j \in \XX$ such that the distributions $\bar{p}_{0, i} = \langle p_{0, ir}, r \in \XX \rangle$ and
$\bar{p}_{0, j} = \langle p_{0, jr}, r \in \XX \rangle$ are singular, i.e., $Q_{0, ij} = \sum_{r \in \XX} \min(p_{0, ir}, p_{0, jr}) = 0$.

Let us consider the random sequence $X^{(N)}_{\e, n} = X_{\e, Nn}, n = 0, 1, \ldots$, for some natural $N > 1$. It is a homogeneous Markov chain, with the initial distribution $\bar{p}$, the phase space
$\XX$, and  the matrix of transition probabilities $\mathbf{P}^N_\e = \| p_{\e, ij}(N) \|$.

 Let us  define the quantities, for $i, j \in \XX$ and $N > 1$,
 \begin{equation}\label{utreda}
Q^{(N)}_{\e, ij} = \sum_{r \in \XX}\min(p_{\e, ir}(N), p_{\e, jr}(N)).
\end{equation}

Let now use the ``multi-step''  coupling algorithm for construction a ``coupling'' Markov chain $Z^{(N)}_{\e, n} =  (X'^{(N)}_{\e, n}, X''^{(N)}_{\e, n}), n = 0, 1, \ldots$, with:

{\bf (i)} the phase space $\ZZ = \XX \times \XX$;

{\bf (ii)} the initial distribution $\bar{P}_{\e} = \langle p_{\e, ij} , (i, j) \in \ZZ \rangle$ constructing according to relation (\ref{solv}), (\ref{solva}), or
(\ref{solvab}) for distributions $\bar{p}' = \bar{p} = \langle p_i, i \in \XX \rangle$ and
$\bar{p}'' = \bar{\pi}_{\e} = \langle \pi_{\e, i}, i \in \XX \rangle$;

{\bf  (iii)} transition probabilities $ P^{(N)}_{\e, ij, rk}$ defined by the following relations, for $(i, j), (r, k) \in \ZZ$:
 \vspace{1mm}

\noindent {\bf (a)} If $Q^{(N)}_{\e, ij} \in (0, 1)$, then,
 \begin{align}\label{gotyrbedgada}
 P_{\e, ij, rk}^{(N)} &  =  \PP \{ X'^{(N)}_{\e, 1} = k, X''^{(N)}_{\e, 1} = r / X'^{(N)}_{\e, 0} = i, X''^{(N)}_{\e, 0} = j \}
\makebox[22mm]{}  \vspace{2mm} \nonumber \\
& = \min(p_{\e, ir}(N), p_{\e, jk}(N)) {\rm I} (r = k) \vspace{2mm} \nonumber \\
& \quad  + \frac{1}{1 - Q^{(N)}_{\e, ij}} ( p_{\e, ir}(N) - \min(p_{\e, ir}(N), p_{\e, jr}(N))) \vspace{2mm} \nonumber \\
& \quad \times ( p_{\e, jk}(N) - \min(p_{\e, ik}(N), p_{\e, jk}(N))),
 \end{align}

\noindent {\bf (b)}  If $Q^{(N)}_{\e, ij} = 1$, then $p_{\e, ir}(N) = p_{\e, jr}(N), r \in \XX$ and,
  \begin{equation}\label{gotyrbedgabda}
 P_{\e, ij, rk}^{(N)} = \min(p_{\e, ir}(N), p_{\e, jk}(N)) {\rm I} (r = k), r, k \in \XX. \makebox[23mm]{}
 \end{equation}

\noindent  {\bf (c)}  If $Q^{(N)}_{\e, ij} = 0$, then $p_{\e, ir}(N) p_{\e, jr}(N) = 0, r \in \XX$ and,
  \begin{equation}\label{gotyrbedgabeda}
  P_{\e, ij, rk}^{(N)}  = p_{\e, ir}(N) p_{\e, jk}(N), r, k \in \XX. \makebox[50mm]{}
 \end{equation}

 The following lemma is the direct corollary of Lemma 7. \vspace{1mm}

{\bf Lemma 8}. {\em $Z^{(N)}_{\e, n} =  (X'^{(N)}_{\e, n}, X''^{(N)}_{\e, n}), n = 0, 1, \ldots$, be a homogeneous Markov chain with the phase space $\ZZ = \XX \times \XX$, the initial distribution $\bar{P_\e}$ and transition probabilities given by relations
{\rm (\ref{gotyrbedgada}) -- (\ref{gotyrbedgabeda})}. Then: \vspace{1mm}

 {\bf (i)} The first component, $X'^{(N)}_{\e, n}, n = 0, 1, \ldots$, is  a homogeneous Markov chain with the phase space $\XX$, the initial distribution $\bar{p}$ and the matrix of transition probabilities $\mathbf{P}^N_\e$.

{\bf (ii)} The second component  $X''^{(N)}_{\e, n}, n = 0, 1, \ldots$  is a homogeneous Markov chain with the phase space $\XX$, the initial distribution $\bar{\pi}_\e$ and the matrix of transition probabilities $\mathbf{P}^N_\e$.

{\bf (iii)} The set $\ZZ_0 = \{ (i, i), i \in \XX \}$ is an absorbing set for the  Markov chain $Z^{(N)}_{\e, n}$, i.e.,  probabilities $P^{(N)}_{\e, ii, rk} = 0$, for $i, r, k  \in \XX, r \neq k$.} \vspace{1mm}

Relation, $\mathbf{A} \mathbf{B} = \mathbf{B}$,  holds for any $m \times m$ stochastic matrix $\mathbf{A} = \| a_{ij} \|$ and $m \times m$  stochastic damping type matrix $\mathbf{B} = \| b_{ij} \|$, with elements $b_{ij} = b_j \geq 0, i, j = 1, \ldots, m$.  Also, matrix
$\mathbf{C} = \mathbf{B} \mathbf{A}$, which has elements, $c_{ij} = c_j = \sum_{k  = 1}^m b_k a_{kj} \geq 0, i, j = 1, \ldots, m$,  is a stochastic damping type matrix, i.e., it has all rows  the same.

Using these remarks, we get the following formula,
\begin{align}\label{formada}
\mathbf{P}^N_\e &  = ((1- \e)\mathbf{P}_0 + \e \mathbf{D})^N
\vspace{1mm} \nonumber \\
& = \mathbf{P}^{N-1}_\e (1- \e)\mathbf{P}_0 + \mathbf{P}^{N-1}_\e  \e \mathbf{D}
\vspace{1mm} \nonumber \\
& =  \mathbf{P}^{N-1}_\e (1- \e)\mathbf{P}_0 + \e \mathbf{D}
\vspace{1mm} \nonumber \\
& = \mathbf{P}^{N-2}_\e (1- \e)^2 \mathbf{P}^2_0 + \mathbf{P}^{N-2}_\e  \e (1- \e) \mathbf{D}\mathbf{P}_0 + \e  \mathbf{D}
\vspace{1mm} \nonumber \\
& = \cdots = (1- \e)^N \mathbf{P}^N_0 + \e  (1- \e)^{N-1}  \mathbf{D}\mathbf{P}^{N-1}_0 + \cdots +  \e  \mathbf{D}.
\end{align}

Using relation (\ref{formada}) and the above remarks we get the following relation,
\begin{align}\label{bedgada}
 Q^{(N)}_\e & = Q(\mathbf{P}^N_\e)  \vspace{1mm} \nonumber \\
 &  \geq (1- \e)^N Q(\mathbf{P}^N_0)  + \e  (1- \e)^{N-1}  Q(\mathbf{D}\mathbf{P}^{N-1}_0) + \cdots + \e  Q(\mathbf{D})
 \vspace{1mm} \nonumber \\
& =   (1- \e)^N Q(\mathbf{P}^N_0)  +  \e  (1- \e)^{N-1}  + \cdots + \e \vspace{1mm} \nonumber \\
& = (1- \e)^N Q(\mathbf{P}^N_0) + 1 - (1- \e)^{N}.
\end{align}

Let us introduce,  for $N \geq 1$, the coefficient of ergodicity,
\begin{equation}\label{jopy}
\Delta_N(\mathbf{P}_0) = (1 - Q(\mathbf{P}^N_0))^{1/N}.
\end{equation}

Note that condition  ${\bf A}_1$ is not required in this theorem.

Below, we count $\Delta_N(\mathbf{P}_0)^0 = 1$, if $\Delta_N(\mathbf{P}_0) = 0$. \vspace{1mm}

{\bf Theorem 6}. {\em The following relation takes place for every $\e \in (0, 1]$ and
$\bar{p} \in {\cal P}_m, j \in \XX, n \geq 0$,}
\begin{equation}\label{bioyt}
 | p_{\e, \bar{p}, j}(n) - \pi_{\e, j} |
\leq ( 1 - Q_{\e, \bar{p}}) \Delta_N(\mathbf{P}_0)^{[n/N]N} (1 - \e)^{[n/N]N}.
\end{equation}

{\bf Proof}. Obviously, for $\bar{p} \in {\cal P}_m, j \in \XX, n \geq 0$,
\begin{equation}\label{statako}
\PP\{ X'^{(N)}_{\e, n} = j \} = p_{\e, \bar{p}, j}(Nn).
\end{equation}

Since, the initial distribution of Markov chain $X''^{(N)}_{\e, n}$ coincides with its stationary distribution, this Markov chain is a stationary random sequence and, thus, for $j \in \XX, n \geq 0$,
\begin{equation}\label{statko}
\PP \{ X''^{(N)}_{\e, n} = j \} = \pi_{\e, j}.
\end{equation}

Let now define a hitting (coupling) time
\begin{equation}\label{statkom}
T^{(N)}_\e = \min(n \geq 0: X'^{(N)}_{\e, n} = X''^{(N)}_{\e, n} ) = \min (n \geq 0: Z^{(N)}_{\e, n} \in \ZZ_0).
\end{equation}

Since $\ZZ_0$ is an absorbing set for the Markov chain $Z_{\e, n}$,  the following relation holds,
\begin{equation}\label{statasko}
\PP \{ Z^{(N)}_{\e, n} \in \ZZ_0,  n \geq T^{(N)}_\e \} = 1.
\end{equation}

Using the above remarks, we get {the following relation, for $j \in \XX, n \geq 0$,
\begin{align}\label{relatdako}
& | p_{\e, \bar{p}, j}(Nn) - \pi_{\e, j} |    = | \PP \{ X'^{(N)}_{\e, n} = j \} -  \PP \{ X''^{(N)}_{\e, n} = j \} |  \vspace{1mm} \nonumber \\
& \quad \quad =   | \PP \{ X'^{(N)}_{\e, n} = j,  X''^{(N)}_{\e, n} \neq j \} -  \PP \{ X'^{(N)}_{\e, n} \neq  j,  X''^{(N)}_{\e, n} = j \} | \vspace{1mm} \nonumber \\
& \quad \quad \leq  \PP \{ X'^{(N)}_{\e, n} = j,  X''^{(2)}_{\e, n} \neq j \} + \PP \{ X'^{(N)}_{\e, n} \neq  j,  X''^{(N)}_{\e, n} = j \}  \vspace{1mm} \nonumber \\
& \quad \quad \leq \PP \{ T^{(N)}_\e > n \}.
\end{align}

Using Lemma 6, we get,  $\bar{p} \in {\cal P}_m, j \in \XX$,
\begin{align}\label{seroft}
| p_{\e, \bar{p}, j}(0) - \pi_{\e, j} |  & \leq  \PP \{ X'^{(N)}_{\e, 0} = j,  X''^{(N)}_{\e, 0} \neq j \} \vspace{1mm} \nonumber \\
& \quad \quad + \PP \{ X'^{(N)}_{\e, 0} \neq  j,  X''^{(N)}_{\e, 0} = j \}
\vspace{1mm} \nonumber \\
&  \leq \PP \{ T^{(N)}_\e  >  0 \} =  1 -  Q_{\e, \bar{p}}.
\end{align}

Also, by continuing the inequalities (\ref{relat}) we get, for  $\bar{p} \in {\cal P}_m, j \in \XX, n \geq 1$,
\begin{align}\label{verwko}
& | p_{\e, \bar{p}, j}(Nn) - \pi_{\e, j} |  \leq  \PP \{ T^{(N)}_\e > n \} \vspace{2mm} \nonumber \\
& \quad \quad =    \sum_{i, j \in \XX} \PP \{  T^{(N)}_\e > n -1,  X'^{(N)}_{\e, n-1} = i, X''^{(N)}_{\e, n-1} = j \}  \vspace{2mm} \nonumber \\
& \quad \quad  \quad  \times  \PP \{ X'^{(N)}_{\e, n} \neq X'^{(N)}_{\e, n} / X'^{(N)}_{\e, n -1} = i,  X'^{(N)}_{\e, n -1} = j \} \vspace{2mm} \nonumber \\
& \quad \quad = \sum_{i, j \in \XX} \PP \{ T^{(N)}_\e > n -1, X'^{(N)}_{\e, n-1} = i, X''^{(N)}_{\e, n-1} = j \} (1 -  Q^{(N)}_{\e, ij})
\vspace{2mm} \nonumber \\
& \quad \quad \leq \PP \{  T^{(N)}_\e > n -1 \} (1 - Q(\mathbf{P}^N_\e)) \vspace{2mm} \nonumber \\
& \quad \quad  \leq  \cdots \leq \PP \{  T_\e > 0 \} (1 - Q(\mathbf{P}^N_\e))^n   \vspace{2mm} \nonumber \\
& \quad \quad   = (1 - Q_{\e, \bar{p}})(1 - Q(\mathbf{P}^N_0))^n (1 - \e)^{Nn} \vspace{2mm} \nonumber \\
& \quad \quad   =  (1 - Q_{\e, \bar{p}})\Delta_N(\mathbf{P}_0))^{Nn} (1 - \e)^{Nn}.
\end{align}

Also, for $\bar{p} \in {\cal P}_m, j \in \XX, n \geq 0$ and $l = 0, \ldots, N -1$,
\begin{align}\label{relatasfa}
& | p_{\e, \bar{p}, j}(Nn+l) - \pi_{\e, j} |   =  |\sum_{k \in \XX} p_{\e, \bar{p}, k}(nN)  p_{\e, kj}(l) - \sum_{k \in \XX} \pi_{\e, k}p_{\e, kj}(l) |  \nonumber \\
& \quad \quad  \quad \leq   \sum_{k \in \XX}   |  p_{\e, \bar{p}, k}(Nn)  - \pi_{\e, k} | p_{\e, kj}(l) \nonumber \\
& \quad \quad  \quad \leq   \max_{k \in \XX}   |  p_{\e, \bar{p}, k}(Nn)  - \pi_{\e, k} |  \nonumber \\
& \quad \quad  \quad \leq (1 - Q_{\e, \bar{p}})\Delta_N(\mathbf{P}_0)^{Nn} (1 - \e)^{Nn}.
\end{align}
\
Inequalities (\ref{seroft}) and (\ref{relatasfa}) imply inequalities given in relation (\ref{bioyt}).
The proof is complete. $\Box$ \vspace{1mm}

The following formula, analogous to (\ref{formad}),  takes place,
\begin{equation}\label{formades}
\Delta_N(\mathbf{P}_0) =  (1 - Q(\mathbf{P}^N_0))^{1/N} =  (\frac{1}{2} \max_{i, j \in \XX} \sum_{k \in \XX} | p_{0, ik}(N) -  p_{0, jk}(N)|)^{1/N}.
\end{equation}

In the case, where condition ${\bf A}_1$ holds, $1 - Q(\mathbf{P}^N_0) \to 0$ as $N \to \infty$,  i.e., the following condition holds for $N$ large enough:
\begin{itemize}
\item [${\bf C}^{(N)}$:] $\Delta_N(\mathbf{P}_0) \in [0, 1)$.
\end{itemize}

Also, condition ${\bf C}^{(N)}$ is, for every $N \geq 1$, sufficient for holding the weaken variant of condition
${\bf A}_1$ mentioned in Remark 4. If the corresponding stationary distribution $\bar{\pi}_0 = \langle \pi_{0, j}, j \in \XX \rangle$ is positive, condition ${\bf C}^{(N)}$ is, for every $N \geq 1$, sufficient for holding condition ${\bf A}_1$.

As in Subsection 7.2, by passing $\e \to 0$ in  the first inequality given in relation (\ref{bioyt})  we get the
inequality holding  for  any  $\bar{p} \in {\cal P}_m, j \in \XX, n \geq 0$,
\begin{equation}\label{hubasqa}
| p_{0, \bar{p}, j}(n) - \pi_{0, j} | \leq
( 1 - Q_{0, \bar{p}}) \Delta_N(\mathbf{P}_0))^{[n/N]N}.
\end{equation}

In the case, where condition  ${\bf A}_1$ holds, the  eigenvalues decomposition representation  (\ref{decom})  implies that the following
relation holds, under  the assumption that the second eigenvalue $\rho_2$ and  the factor which stays in front of $|\rho_2 |$ in relation (\ref{formadesam}) are not equal $0$,
\begin{align}\label{formadesam}
 \Delta_N \big( \mathbf{P}_0) &  \sim (\frac{1}{2} \max_{i, j \in \XX}
\sum_{k \in \XX} | \sum_{2 \leq l \leq m: |\rho_l| = |\rho_2|} ( \pi_{0, ik}[l] - \pi_{0, jk}[l]) | \big)^{1/N} |\rho_2| \vspace{2mm} \nonumber \\
&  \to  |\rho_2| \ {\rm as} \ N \to \infty.
\end{align}

Relations (\ref{hubasqa}) and (\ref{formadesam}) show that, in the case where condition ${\bf A}_1$ holds,   the  ``coupling''  upper bounds for rate of convergence in individual ergodic relations given in Theorem 6 are asymptotically equivalent with analogous upper bounds, which can be obtained with the use of eigenvalue decomposition  representation for transition probabilities. At the same time, computing of coefficients of ergodicity $\Delta_N(\mathbf{P}_0)$ does not require solving of the polynomial
equation, $\det(\rho \mathbf{I} - \mathbf{P}_0) = 0$, that is required for finding eigenvalues. \vspace{1mm}

 {\bf 7.4 Coupling  for singularly perturbed Markov chains  with damping component}. In this subsection, we present  coupling algorithms and get the effective upper bound for the rate of convergence in ergodic theorems for singularly perturbed Markov chains with damping component. Thus, we assume that following weaken variant of condition  ${\bf B}_1$  holds:
 \begin{itemize}
\item [${\bf B}_4$:] The phase space $\XX = \cup_{j = 1}^h \XX^{(j)}$, where: (a) $X^{(j)}, j = 1, \ldots, h$ are non-intersecting subsets of $\XX$, (b)
$X^{(j)}, j = 1, \ldots, h$ are non-empty, closed classes of  states for the Markov chain $X_{0, n}$.
\end{itemize}

Let us introduce the discrete distributions $\bar{f}_{\bar{p}} = \langle f^{(j)}_{\bar{p}}, j = 1, \ldots, h \rangle$,  where $f^{(j)}_{\bar{p}} = \sum_{k \in \XX^{(j)}} p_k,  , j = 1, \ldots, h$
and $\bar{p}^{(j)} = \langle p^{(j)}_{k} =  p_k / f^{(j)}_{\bar{p}},  k \in \XX^{(j)} \rangle$, for $j \in \HH_{\bar{p}} = \{j: 1 \leq j \leq h,   f^{(j)}_{\bar{p}}  > 0 \}$.

We also use the above notations for the case, where distribution $\bar{p}$ coincides with distributions $\bar{d} = \langle d_k, k \in \XX �\rangle$,  $\bar{\pi}_\e = \langle \pi_{\e, k}, k \in \XX �\rangle$, and $\bar{\pi}_{0, \bar{p}} = \langle \pi_{0, \bar{p},  k}, k \in \XX �\rangle$.

Let $\e \in (0, 1]$, and $\bar{\pi}_\e = \langle \pi_{\e, k}, k \in \XX \rangle$ be the stationary distribution for the Markov chain $X_{\e, n}$.

The following representation takes place,  for  $k \in \XX^{(j)}, j = 1, \ldots, h$,
\begin{equation}\label{reprad}
\pi_{\e, k}  = \e \sum_{l = 0}^{\infty}  p_{0, \bar{d}, k}(l)(1- \e)^l  = f^{(j)}_{\bar{d}} \pi^{(j)}_{\e, k},
\end{equation}
where, for $k \in \XX^{(j)}, j = 1, \ldots, h$,
\begin{equation}\label{hopwe}
\pi^{(j)}_{\e, k}  =   \e \sum_{l = 0}^{\infty}  p_{0, \bar{d}^{(j)}, k}(l)(1- \e)^l.
\end{equation}

It is readily seen that $\bar{\pi}^{(j)}_\e = \langle \pi^{(j)}_{\e, k}, k \in \XX \rangle$ is, for every $j = 1, \ldots, h$, the stationary distribution for the Markov chain $X^{(j)}_{\e, n}$, with the  phase space $\XX^{(j)}$ and the matrix of transition probabilities $\mathbf{P}_{\e, j} =
(1- \e) \mathbf{P}_{0, j} + \e \mathbf{D}_j$, where
$\mathbf{P}_{0, j} = \| p_{0, rk} \|_{r, k \in \XX^{(j)}}$ is, according to condition ${\bf B}_4$, a stochastic matrix, while $\mathbf{D}_j = \| d^{(j)}_{r k} \|_{r, k \in \XX^{(j)}}$ is the stochastic damping matrix
with elements $d^{(j)}_{rk} =  d^{(j)}_{k} = d_k / f^{(j)}_{\bar{d}}, k, r \in \XX^{(j)}$.

Note that relation (\ref{reprad}) implies that, for $j = 1, \ldots, h$,
\begin{equation}\label{reprado}
f^{(j)}_{\bar{\pi}_\e} = \sum_{k \in  \XX^{(j)}} \pi_{\e, k} = \sum_{k \in  \XX^{(j)}} f^{(j)}_{\bar{d}} \pi^{(j)}_{\e, k} = f^{(j)}_{\bar{d}}.
\end{equation}

Let us also introduce quantities, for $j = 1, \ldots, h, N \geq 1$,
\begin{equation}\label{initacaba}
Q^{(j)}_N(\mathbf{P}_0) = \min_{k, r \in \XX^{(j)}} \sum_{i \in \XX^{(j)}} p_{0, ki}(N) \wedge p_{0, ri}(N),
\end{equation}
and
\begin{equation}\label{initacabanu}
\Delta^{(j)}_N(\mathbf{P}_0) = (1 - Q^{(j)}_N(\mathbf{P}_0))^{1/N}.
\end{equation}

It is useful to note that  $\Delta^{(j)}_N(\mathbf{P}_0) \to 0$ as $N \to \infty$, and, thus,  $\Delta^{(j)}_N(\mathbf{P}_0) < 1$, for all $N$ large enough, if condition ${\bf A}_1$ holds for the Markov chain $X^{(j)}_{0, n}$.

Also the weaken variant of condition ${\bf A}_1$, pointed out in Remark 4,  holds for the Markov chain $X^{(j)}_{0, n}$, if $\Delta^{(j)}_N(\mathbf{P}_0) < 1$, for some $N \geq 1$.

In this case, the Markov chain $X^{(j)}_{0, n}$ has the stationary distribution $\bar{\pi}^{(j)}_{0} =
\langle \pi^{(j)}_{0, k}, k \in \XX^{(j)} \rangle$.

Let us now assume that the following condition holds, for some $N \geq 1$:
\begin{itemize}
\item [${\bf D}^{(N)}$:] {\bf (a)} Condition ${\bf B}_4$ holds, {\bf (b)} $\Delta^{(j)}_{N}(\mathbf{P}_0) < 1$, for $j = 1, \ldots, h$.
\end{itemize}

According to Remark 4,  the Markov chain $X_{0, n}$ has the stationary distribution $\bar{\pi}_{0, \bar{p}} = \langle \pi_{0, \bar{p}, k}$, $k \in \XX^{(j)} \rangle$, for which the following relation holds, for $k \in \XX^{(j)}, j = 1, \ldots, h$,
\begin{equation}\label{hotyre}
\pi_{0, \bar{p}, k} = f^{(j)}_{\bar{p}} \pi^{(j)}_{0, k}.
\end{equation}

It is useful to note that transition probabilities $p_{0, \bar{p}, k}(n) = 0, n \geq 0$ and stationary probabilities  $\pi_{0, \bar{p}, k} = 0$, for
$k \in \XX^{(j)}, j \notin \HH_{\bar{p}}$.

Moreover,  if $k \in \XX^{(j)}$, for some $j = 1, \ldots, h$, then,
\begin{equation}\label{kuk}
p_{0, \bar{p}, k}(n) = f^{(j)}_{\bar{p}} \PP_{\bar{p}^{(j)}} \{ X_{0, n}^{(j)} = k \}.
\end{equation}

Below, we count $\Delta^{(j)}_N(\mathbf{P}_0)^0 = 1$, if $\Delta^{(j)}_N(\mathbf{P}_0) = 0$. \vspace{1mm}

{\bf Theorem  7}. {\em Let condition ${\bf D}^{(N)}$  holds. Then,  if $\bar{p} \in {\cal P}_m$, and $k \in \XX^{(j)}$, for some $j = 1, \ldots, h$, the following propositions take place, for every $\e \in (0, 1]$ and $n \geq 0$}: \vspace{1mm}
\begin{align}\label{husner}
 | p_{\e, \bar{p}, k}(n) - \pi_{\e,  k} |  & \leq \big( ( f^{(j)}_{\bar{d}}( 1 - Q(\bar{\pi}_\e^{(j)}, \bar{\pi}_0^{(j)}))
\vspace{3mm}  \nonumber \\
& \quad  + f^{(j)}_{\bar{p}}( 1 - Q(\bar{p}^{(j)}, \bar{\pi}_0^{(j)}))) \Delta^{(j)}_N(\mathbf{P}_0)^{[n/N]N}
\makebox[15mm]{}\vspace{3mm}  \nonumber \\
& \quad + | f^{(j)}_{\bar{p}} - f^{(j)}_{\bar{d}}| \pi^{(j)}_{0, k}  \big)(1 - \e)^{n}.
\end{align}

{\bf Proof}. By using the renewal  equation (\ref{gotyrtbaraka}) and the renewal type type relation (\ref{gotyrtbarakaba}) and taking into account stationarity of the Markov chain $X_{\e, n}$, with the initial distribution $\bar{\pi}_\e$, condition ${\bf B}_4$,  and
relations  (\ref{reprad}), (\ref{hotyre}) and (\ref{kuk}), we get the following relation, for $k \in \XX^{(j)}$ and $n \geq 0$,
\begin{align}\label{bopkur}
| p_{\e, \bar{p}, k}(n) - \pi_{\e, k} | & = | p_{\e, \bar{p}, k}(n) - p_{\e, \bar{\pi}_\e, k}(n) |  \vspace{2mm} \nonumber \\
& = | p_{0, \bar{p}, k}(n) - p_{0, \bar{\pi}_\e, k}(n) | (1 - \e)^n \vspace{2mm} \nonumber \\
& =  | f^{(j)}_{\bar{p}} \PP_{\bar{p}^{(j)}} \{ X_{0, n}^{(j)} = k \} - f^{(j)}_{\bar{d}} \PP_{\bar{\pi}^{(j)}_{\e}} \{ X_{0, n}^{(j)} = k \} | (1 - \e)^{n} \makebox[2mm]{} \vspace{2mm} \nonumber \\
&  \leq \big(f^{(j)}_{\bar{p}}  | \PP_{\bar{p}^{(j)}} \{ X_{0, n}^{(j)} = k \} -   \pi^{(j)}_{0, k} |  \vspace{2mm} \nonumber \\
& \quad +   f^{(j)}_{\bar{d}}| \PP_{\bar{\pi}^{(j)}_{\e}} \{ X_{0, n}^{(j)} = k \}  -   \pi^{(j)}_{0, k} | \vspace{2mm} \nonumber \\
&  \quad +  | f^{(j)}_{\bar{p}} - f^{(j)}_{\bar{d}} | \pi^{(j)}_{0, k} \big) (1 - \e)^{n}.
\end{align}

Condition ${\bf D}^{(N)}$ let one apply inequality (\ref{hubasqa}) to the Markov chain $X_{0, n}^{(j)}$ and transforms  relation (\ref{bopkur}) to the form of relation (\ref{husner}). $\Box$ \vspace{1mm}

As in Subsection 7.3,  by passing $\e \to 0$ in  the first inequality given in relation (\ref{husner}) and taking into account Lemma 5
(which implies that  $Q(\bar{\pi}_\e^{(j)}, \bar{\pi}_0^{(j)}) \to 1$ as $\e \to 0$),   we get the following relation holding  for  any $\bar{p} \in {\cal P}_m, k \in \XX^{(j)}, j = 1, \ldots, h, n \geq 0$,
\begin{align}\label{husnerve}
& | p_{0, \bar{p}, k}(n) - \pi_{0, \bar{d},  k} |
\vspace{2mm} \nonumber \\
& \quad \quad \leq f^{(j)}_{\bar{p}}( 1 - Q(\bar{p}^{(j)}, \bar{\pi}_0^{(j)})) \Delta^{(j)}_N(\mathbf{P}_0)^{[n/N]N}
 + | f^{(j)}_{\bar{p}} - f^{(j)}_{\bar{d}}| \pi^{(j)}_{0, k}.
\end{align}

This relation let one extend the conclusion made in Remark 2.

If $f^{(j)}_{\bar{p}} = f^{(j)}_{\bar{d}}$, then, for $k \in \XX^{(j)}$, probabilities $p_{0, \bar{p}, k}(n) \to  \pi_{0, \bar{d},  k}$ as $n \to \infty$, with the exponential rate.

If $f^{(j)}_{\bar{p}} \neq f^{(j)}_{\bar{d}}$, then,  for $k \in \XX^{(j)}$ (as follows from Lemma 5 and relations  (\ref{reprad}) and (\ref{hotyre})),  probabilities $p_{0, \bar{p}, k}(n) \to  \pi_{0, \bar{p},  k} =  \pi_{0, \bar{d},  k} +  (f^{(j)}_{\bar{p}} - f^{(j)}_{\bar{d}}) \pi^{(j)}_{0, k}$ as $n \to \infty$, with the exponential rate. More precisely, relations (\ref{hubasqa}), (\ref{hotyre}),  and (\ref{kuk}) imply that the following relation holds, for $k \in \XX^{(j)}, n \geq 0$,
\begin{align}\label{husnervem}
| p_{0, \bar{p}, k}(n) - \pi_{0, \bar{p},  k} | & =
| f^{(j)}_{\bar{p}} \PP_{\bar{p}^{(j)}} \{ X_{0, n}^{(j)} = k \} -  f^{(j)}_{\bar{p}} \pi^{(j)}_{0, k}| \vspace{2mm} \nonumber \\
& \leq f^{(j)}_{\bar{p}}( 1 - Q(\bar{p}^{(j)}, \bar{\pi}_0^{(j)})) \Delta^{(j)}_N(\mathbf{P}_0)^{[n/N]N}.
\end{align}
\vspace{1mm}

{\bf 8. Ergodic theorems for perturbed Markov chains with  \\ \makebox[11mm]{}  damping component in the triangular array mode} \\

In this section, we present ergodic theorems for Markov chains with damping component in the triangular array mode, where time $n$ tends to infinity and damping parameter $\e$ tends to zero simultaneously. In this mode, the asymptotic behaviour of transition probabilities  \vspace{1mm}

{\bf 8.1 Asymptotic behaviour of transition probabilities  in the triangular array mode}.  In this mode, the asymptotic behaviour
of transition  probabilities $p_{\e, \bar{p}, k}(n)$ is studied when  time $n \to \infty$ and the damping parameter  $\e \to 0$,  simultaneously.
One can assume that time $n = n_\e$ depends on $\e$ in such way that $n_\e \to \infty$ as $\e \to 0$.

Asymptotic behaviour of transition probabilities $p_{\e, \bar{p}, k}(n_\e)$ as $\e \to 0$ should be compared with repeated limits for
$p_{\e, \bar{p}, k}(n)$ as $n \to \infty$ and, then, $\e \to 0$ and, vice versa,  as $ \e \to 0$ and, then $n \to \infty$.

Let as first  consider the case, where condition ${\bf A}_1$ holds.

In this case, according to Lemma 2,  transition probabilities $p_{\e, \bar{p}, k}(n) \to \pi_{\e, k}$ as $n \to \infty$, for every $\e \in (0, 1]$. Also, according to Lemma 4, stationary probabilities $\pi_{\e, k}  \to \pi_{0, k}$ as $\e \to 0$, where
$\pi_{0, k}$  are stationary probabilities for the Markov chain $X_{0, n}$.

As follows from relation (\ref{transava}), in this case, transition probabilities $p_{\e, \bar{p}, k}(n)$ $\to p_{0, \bar{p}, k}(n)$ as $\e \to 0$, for any $n \geq 0$. Also, according to relation (\ref{erg}), transition probabilities $p_{0, \bar{p}, k}(n) \to \pi_{0, k}$ as $n \to 0$.

The coincidence, of the repeated limits for transition probabilities $p_{\e, \bar{p}, k}(n)$ let one expect that, in this case,  transition  probabilities $p_{\e, \bar{p}, k}(n_\e)$ should converge to  stationary probabilities $\pi_{0, k}$  as $\e \to 0$, for an arbitrary $n_\e \to \infty$ as $\e \to 0$.

The situation is different, in the case,  where condition ${\bf B}_1$ holds.

Again, Lemma 2 implies that  transition probabilities $p_{\e, \bar{p}, k}(n)$ $\to \pi_{\e, k}$ as $n \to \infty$, for every $\e \in (0, 1]$. Also,  according to Lemma 4, stationary probabilities $\pi_{\e, k}  \to  \pi_{0, \bar{d},  k}$ as $\e \to 0$, where
$\pi_{0, \bar{d},  k}$  are stationary probabilities of the Markov chain $X_{0, n}$, with the initial distribution $\bar{d}$.

Again, according to relation (\ref{transava}), transition probabilities $p_{\e, \bar{p}, k}(n) \to p_{0, \bar{p}, k}(n)$ as $\e \to 0$, for any $n \geq 0$. However, as Lemma 3 implies that, in this case, transition probabilities $p_{0, \bar{p}, k}(n) \to \pi_{0, \bar{p}, k}$ as $n \to \infty$.

Thus, the repeated limits for transition probabilities $p_{\e, \bar{p}, k}(n)$ coincides, under additional assumption that the initial distribution
$\bar{p} = \bar{d}$. But, they may not coincide if $\bar{p} \neq \bar{d}$.

The transition probabilities  $p_{\e, \bar{p}, k}(n)$ satisfy the renewal type relation (\ref{gotyrtbarakaba}). This relation shows that
the initial state $X_{\e, 0}$ influences the behaviour of the Markov chain $X_{\e, n}$  only up to the first damping regeneration time $T_{\e, 1}$.

The random variables $T_{\e, 1}$ has the geometric distribution with parameter $\e$. Its  rate of growth is $\e^{-1}$ as $\e \to 0$. Moreover, random variables $\e T_{\e.1} \stackrel{d}{\longrightarrow} \nu(1)$ as $\e \to 0$, where $\nu(1)$ is a random variable,  which has exponential distribution with parameter $1$.

This hints one that it would be natural to study the asymptotic behaviour of transition  probabilities $p_{\e, \bar{p}, k}(n_\e)$
for $n_\e \to \infty$ as $\e \to 0$ such that $\e n_\e \to t \in [0, \infty]$ as $\e \to 0$.

Moreover, it can be expected that in two extremal cases, where $t = \infty$ or $t = 0$, the transition probabilities
$p_{\e, \bar{p}, k}(n_\e)$ converge as $\e \to 0$ to the corresponding repeated limits, respectively,
$\pi_{0, \bar{d},  k}$ or  $\pi_{0, \bar{p},  k}$.

The question arises about the asymptotic behaviour of  transition probabilities $p_{\e, \bar{p}, k}(n_\e)$ in the intermediate case,  where
the above limit $t  \in (0, \infty)$.

In conclusion of this informal discussion, we would like to mention works \cite{En2, EnS, GySi5, GS4, Pe2, Si1, Si2, Si8, SiP},  which contain
results related to ergodic theorems in triangular array mode and to so-called quasi-stationary ergodic theorems  for perturbed regenerative processes, Markov chains and semi-Markov processes.  

\vspace{1mm}

{\bf 8.2 Ergodic theorems for regularly perturbed Markov chains with damping component in the triangular array mode}. The following
theorem takes place. \vspace{1mm}

{\bf Theorem 8}. {\em Let condition  ${\bf C}^{(N)}$ holds for some $N \geq 1$. Then, for any $n_\e \to \infty$ as $\e \to 0$ and $\bar{p} \in {\cal P}_m, k \in \XX$,}
\begin{equation}\label{quasir}
p_{\e, \bar{p}, k}(n_\e) \to \pi_{0, k} \ {\rm as} \ \e \to 0.
\end{equation}

{\bf Proof}. Using the renewal type relation  (\ref{gotyrtbarakaba}) and  inequality  (\ref{hubasqa}), we get the following relation holding for any $n_\e \to \infty$ as $\e \to 0$ and
$k \in \XX$,
\begin{align}\label{gotyrabamol}
|p_{\e, \bar{p}, k}(n_\e)  - \pi_{0, k} |  & = |p_{\e, \bar{p}, k}(n_\e)  - \pi_{0, k} (1- \e)^{n_\e} -  \pi_{0, k}(1-  (1 - \e)^{n_\e}) |
\vspace{3mm} \nonumber \\
&  =   | (p_{0, \bar{p}, k}(n_\e) -\pi_{0, k}) (1- \e)^{n_\e}  \vspace{2mm} \nonumber \\
& \quad  +  \sum_{l = 1}^{n_\e}  ( p_{0, \bar{d}, k}(n_\e -l) -  \pi_{0, k}) \e (1- \e)^{l - 1}| \vspace{2mm} \nonumber \\
&  \leq ( 1 - Q(\bar{p}, \bar{\pi}_0))  \Delta_N(\mathbf{P}_0)^{[n_\e/N]N} \vspace{2mm} \nonumber \\
& \quad + ( 1 - Q(\bar{d}, \bar{\pi}_0))  \e \sum_{l = 1}^{n_\e} \Delta_N(\mathbf{P}_0)^{[(n_\e - l)/N]N} \vspace{2mm} \nonumber \\
&  \leq ( 1 - Q(\bar{p}, \bar{\pi}_0))   \Delta_N(\mathbf{P}_0)^{[n_\e/N]N}  \vspace{2mm} \nonumber \\
& \quad  + ( 1 - Q(\bar{d}, \bar{\pi}_0)) \e N(1 -  \Delta_N(\mathbf{P}_0)^N)^{-1}  \vspace{2mm} \nonumber \\
& \to 0 \ {\rm as} \ \e \to 0.
\end{align}

The proof is complete. $\Box$ \vspace{1mm}

 {\bf Remark 6}. Relation {\rm (\ref{gotyrabamol})} gives, in fact, explicit upper bounds for the rate of convergence in ergodic relation given in Theorem {\rm 8}. \vspace{1mm}

{\bf 8.3 Ergodic theorems for singularly perturbed Markov chains with damping component in the triangular array mode}.
As follows from remarks made in Subsection 8.1, the ergodic relations take in  this case much more complex forms. \vspace{1mm}

{\bf Theorem 9}. {\em Let condition  ${\bf D}^{(N)}$ holds for some $N \geq 1$. Then the following propositions take place{\rm :}  \vspace{1mm}

{\bf (i)} If $n_\e \to \infty$ and $\e n_\e \to \infty$ as $\e \to 0$, then, for $\bar{p} \in {\cal P}_m, k \in \XX_j, j = 1, \ldots, h$,
\begin{equation}\label{first}
p_{\e, \bar{p}, k}(n_\e) \to \pi_{0, \bar{p}, k}(\infty) = \pi_{0, \bar{d}, k} \ {\rm as} \ \e \to 0. \makebox[32mm]{}
\end{equation}

{\bf (ii)} If $n_\e \to \infty$ and  $\e n_\e \to t \in  (0, \infty)$ as $\e \to 0$, then, for $\bar{p} \in {\cal P}_m, k \in \XX_j, j = 1, \ldots, h$,
\begin{equation}\label{firsta}
p_{\e, \bar{p}, k}(n_\e) \to \pi_{0, \bar{p}, k}(t) = \pi_{0, \bar{p}, k} e^{- t}  + \pi_{0, \bar{d},  k}(1 - e^{- t})  \ {\rm as} \ \e \to 0.
\end{equation}

{\bf (iii)} If $n_\e \to \infty$ and $\e n_\e \to 0$ as $\e \to 0$, then, for $\bar{p} \in {\cal P}_m, k \in \XX_j, j = 1, \ldots, h$,
\begin{equation}\label{firstab}
p_{\e, \bar{p}, k}(n_\e) \to \pi_{0, \bar{p}, k}(0) = \pi_{0, \bar{p}, k}   \ {\rm as} \ \e \to 0. \makebox[36mm]{}
\end{equation}
}
{\bf Proof}. The renewal type relation  (\ref{gotyrtbarakaba}) written for $n = n_\e$ takes the following form, for $\bar{p} \in {\cal P}_m, k \in \XX^{(j)}, k = 1, \ldots, h$,
\begin{equation}\label{gotyraba}
p_{\e, \bar{p}, k}(n_\e) = p_{0, \bar{p}, k}(n_\e)(1- \e)^{n_\e}   + \sum_{l = 1}^{n_\e}  p_{0, \bar{d}, k}(n_\e - l) \e (1- \e)^{l - 1}.
\end{equation}

By applying inequality (\ref{husnervem}) to the transition probabilities appearing in  the  above renewal type relation, we get the following inequality,
\begin{align}\label{gotyrabamo}
& |p_{\e, \bar{p}, k}(n_\e)  - \pi_{0, \bar{p}, k} (1- \e)^{n_\e} -  \pi_{0, \bar{d}, k}(1-  (1 - \e)^{n_\e}) |
\vspace{3mm} \nonumber \\
& \quad \quad =   | (p_{0, \bar{p}, k}(n_\e) -\pi_{0, \bar{p}, k}) (1- \e)^{n_\e}  \vspace{2mm} \nonumber \\
& \quad \quad \quad +  \sum_{l = 1}^{n_\e}  ( p_{0, \bar{d}, k}(n_\e -l) -  \pi_{0, \bar{d}, k}) \e (1- \e)^{l - 1}| \vspace{3mm} \nonumber \\
& \quad \quad \leq f^{(j)}_{\bar{p}}( 1 - Q(\bar{p}^{(j)}, \bar{\pi}_0^{(j)}))  \Delta^{(j)}_N(\mathbf{P}_0)^{[n_\e/N]N} +
  \vspace{2mm} \nonumber \\
& \quad \quad \quad +  f^{(j)}_{\bar{d}}( 1 - Q(\bar{d}^{(j)}, \bar{\pi}_0^{(j)}))\e \sum_{l = 1}^{n_\e} \Delta^{(j)}_N(\mathbf{P}_0)^{[(n_\e - l)/N]N}  \vspace{2mm} \nonumber \\
& \quad \quad \leq f^{(j)}_{\bar{p}}( 1 - Q(\bar{p}^{(j)}, \bar{\pi}_0^{(j)}))  \Delta^{(j)}_N(\mathbf{P}_0)^{[n_\e/N]N}
\vspace{2mm} \nonumber \\
& \quad \quad \quad  + f^{(j)}_{\bar{d}}( 1 - Q(\bar{d}^{(j)}, \bar{\pi}_0^{(j)})) \e N(1 -  \Delta^{(j)}_N(\mathbf{P}_0)^N)^{-1}.
\end{align}

Let us introduce function $R_\e(t) = | (1 - \e)^{n_\e} - e^{-t}|, t \in [0, \infty]$.

If $n_\e \to \infty$ and  $\e n_\e \to t \in  [0, \infty]$ as $\e \to 0$, then,
\begin{equation}\label{refa}
R_\e(t) \to 0 \ {\rm as} \ \e \to 0.
\end{equation}

The following inequality takes place,
\begin{align}\label{gobamo}
& | \pi_{0, \bar{p}, k} e^{-t} + \pi_{0, \bar{d}, k} (1 -  e^{-t}) -  \pi_{0, \bar{p}, k} (1- \e)^{n_\e} -  \pi_{0, \bar{d}, k}(1-  (1 - \e)^{n_\e}) | \vspace{2mm} \nonumber \\
& \quad \quad \leq |\pi_{0, \bar{p}, k} - \pi_{0, \bar{d}, k}| R_\e(t)  \vspace{2mm} \nonumber \\
& \quad \quad =  |f^{(j)}_{\bar{p}} - f^{(j)}_{\bar{d}} | \pi^{(j)}_{0, k} R_\e(t).
 \end{align}

Relations  (\ref{gotyrabamo}) and (\ref{gobamo}) obviously imply that the following relation holds, for $k \in \XX^{(j)}, k = 1, \ldots, h$,  if  $n_\e \to \infty$ and  $\e n_\e \to t \in  [0, \infty]$ as $\e \to 0$,

\begin{align}\label{gotyrao}
& |p_{\e, \bar{p}, k}(n_\e)  - \pi_{0, \bar{p}, k}(t) | \vspace{2mm} \nonumber \\
& \quad \quad \leq f^{(j)}_{\bar{p}}( 1 - Q(\bar{p}^{(j)}, \bar{\pi}_0^{(j)}))  \Delta^{(j)}_N(\mathbf{P}_0)^{[n_\e/N]N}  \vspace{2mm} \nonumber \\
& \quad \quad \quad
+ f^{(j)}_{\bar{d}}( 1 - Q(\bar{d}^{(j)}, \bar{\pi}_0^{(j)})) \e N(1 -  \Delta^{(j)}_N(\mathbf{P}_0)^N)^{-1}   \vspace{2mm} \nonumber \\
& \quad \quad \quad +  |f^{(j)}_{\bar{p}} - f^{(j)}_{\bar{d}} | \pi^{(j)}_{0, k} R_\e(t)  \to 0 \ {\rm as} \ \e \to 0.
\end{align}

This relation proves the theorem. $\Box$ \vspace{1mm}

{\bf Remark 7}. Inequality {\rm (\ref{gotyrao})}  gives, in fact, explicit upper bounds for the rate of convergence in ergodic relation given in Theorem {\rm 9}.  Of course, it is possible to get some simple explicit upper bounds for  $R_\e(t)$ in terms of  quantities $\e n_\e$ and $t$. \\

{\bf 9. Examples and results of numerical experiments} \\

In this section we present examples illustrating the theoretical results of the present paper using some related numerical experiments. The first example demonstrates results associated to regularly perturbed Markov chains. For instance, asymptotic representation of stationary distribution, estimation of coefficient of ergodicity and associated upper bounds for the rate of convergence for stationary distribution. The second example focuses on the advantage of higher order terms of asymptotic expansion in regularly perturbed Markov chains. Lastly, we present information network that capture results associated with singularly perturbed Markov chains with damping component in the triangular array mode.  \vspace{1mm}

{\bf 9.1 A simple information network with a complete link graph}. In this subsection we consider the Markov chain associated
with a simple information network with five nodes and a complete link graph. We also restrict consideration  by the model with the simplest
damping matrix with all equal elements. The matrices $\mathbf{P}_0$ and $\mathbf{D}$ have the following forms,
\begin{equation}\label{amanomop}
\mathbf{P}_0 = \left\|
\begin{array}{ccccccc}
\frac{1}{5} & \frac{1}{5} & \frac{1}{5} & \frac{1}{5} & \frac{1}{5} \vspace{1mm}  \\
\frac{1}{4} & 0 & \frac{1}{4}  & \frac{1}{4} & \frac{1}{4}  \vspace{1mm} \\
0 & \frac{1}{3} & 0 & \frac{1}{3} &  \frac{1}{3} \vspace{1mm} \\
0 & \frac{1}{3} & \frac{1}{3} & 0 & \frac{1}{3}  \vspace{1mm} \\
0 & \frac{1}{3} & \frac{1}{3} & \frac{1}{3} & 0   \\
\end{array}
\right\|,  \ \mathbf{D}   = \left\|
\begin{array}{ccccccc}
\frac{1}{5} & \frac{1}{5} & \frac{1}{5} & \frac{1}{5} & \frac{1}{5} \vspace{1mm}  \\
\frac{1}{5} & \frac{1}{5} & \frac{1}{5} & \frac{1}{5} & \frac{1}{5}  \vspace{1mm} \\
\frac{1}{5} & \frac{1}{5} & \frac{1}{5} & \frac{1}{5} & \frac{1}{5}  \vspace{1mm} \\
\frac{1}{5} & \frac{1}{5} & \frac{1}{5} & \frac{1}{5} & \frac{1}{5} \vspace{1mm}  \\
\frac{1}{5} & \frac{1}{5} & \frac{1}{5} & \frac{1}{5} & \frac{1}{5}   \\
\end{array}
\right\|.
\end{equation}

\vspace{1mm}

The corresponding link graph is presented in Figure 1.


\begin{figure}[H]
\begin{center}
\includegraphics [width=.5\textwidth]{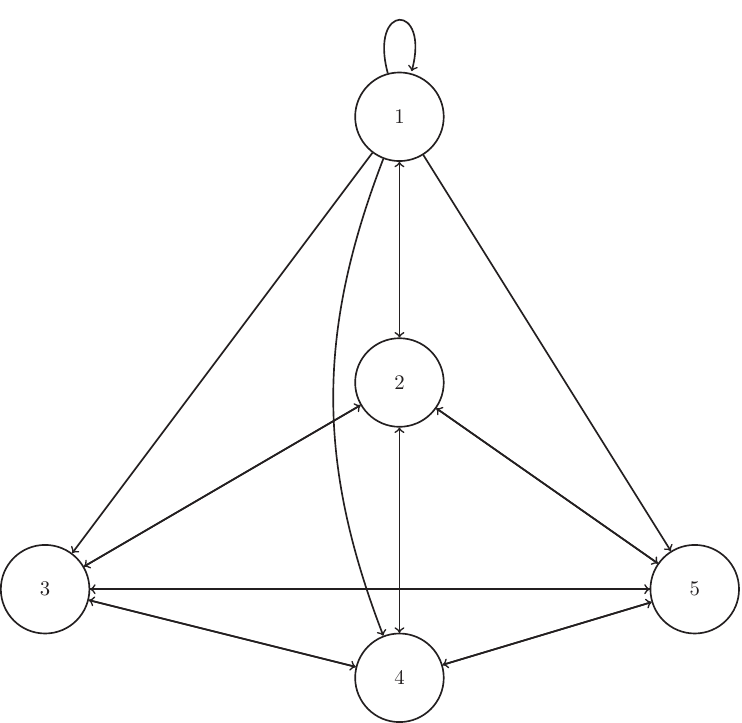}
%
%
%
%
%
\end{center}
\caption{A links graph of simple information network} \label{fig:2}
\end{figure}

In this case, condition $\mathbf{A}_1$ holds.  

The eigenvalues of matrix $\mathbf{P}_0$ (computed by solving equation,  $\det(\rho \mathbf{I} - \mathbf{P}_0)$ $= 0$), are,
\begin{equation}\label{eigen}
 \rho_1 = 1, \  \rho_2 = \rho_3 = - \frac{1}{3}, \  \rho_4 = - \frac{1}{15} - \frac{\sqrt{34}}{30},  \  \rho_5 = - \frac{1}{15} + \frac{\sqrt{34}}{30}.
\end{equation}

The stationary projector matrix $\mathbf{\Pi} = \| \pi_{0, ij} \|$, with  elements $\pi_{0, ij} = \pi_{0, j}, i, j = 1, \ldots, m$ (computed by solving the system of linear equations (\ref{ergsys})), and matrices   $\bar{\mathbf{\Pi}}_k = \| \bar{\pi}_{ij}[k] \|, k = 2, 3, 4$ (computed with the use of the algorithm described in \cite{Fel}), take the following forms,
{\footnotesize
\begin{equation}\label{pitom}
 \mathbf{\Pi} = \left\|
\begin{array}{rrrrrr}
\frac{5}{66} & \frac{8}{33} & \frac{5}{22} & \frac{5}{22} & \frac{5}{22} \vspace{1mm} \\
\frac{5}{66} & \frac{8}{33} & \frac{5}{22} & \frac{5}{22} & \frac{5}{22} \vspace{1mm} \\
\frac{5}{66} & \frac{8}{33} &  \frac{5}{22} & \frac{5}{22} & \frac{5}{22} \vspace{1mm} \\
\frac{5}{66} & \frac{8}{33} &  \frac{5}{22} & \frac{5}{22} & \frac{5}{22} \vspace{1mm} \\
\frac{5}{66} & \frac{8}{33} &  \frac{5}{22} & \frac{5}{22} & \frac{5}{22}
\end{array}
\right\|, \, \bar{\mathbf{\Pi}}_{2}= \left\|
 \begin{array}{rrrrrr}
0 & 0 & 0 & 0 & 0 \vspace{1mm} \\
0 & 0 & 0 & 0 & 0 \vspace{1mm} \\
0 & 0 & \frac{2}{3} & -  \frac{1}{3} &  -  \frac{1}{3}  \vspace{1mm} \\
0 & 0 & - \frac{1}{3} &   \frac{2}{3} &  -  \frac{1}{3} \vspace{1mm} \\
0 & 0 & -  \frac{1}{3} & -  \frac{1}{3} &  \frac{2}{3}
 \end{array}
\right\|,
\end{equation}}
and
{\footnotesize
\begin{equation}\label{pitoma}
 \bar{\mathbf{\Pi}}_{3} = \left\|
\begin{array}{rrrrrr}
-a_{1} & -b_{1} & c_{1} & c_{1} & c_{1} \vspace{1mm} \\
a_{2} & b_{2} & -c_{2} & -c_{2} & -c_{2} \vspace{1mm} \\
-a_{3} & -b_{3} & c_{3} & c_{3} & c_{3} \vspace{1mm} \\
-a_{3} & -b_{3} & c_{3} & c_{3} & c_{3} \vspace{1mm} \\
-a_{3} & -b_{3} & c_{3} & c_{3} & c_{3}
\end{array}
\right\|,  \bar{\mathbf{\Pi}}_{4} = \left\|
\begin{array}{rrrrrr}
e_{1} & f_{1} & -g_{1} & -g_{1} & -g_{1} \vspace{1mm} \\
e_{2} & -f_{2} & g_{2} & g_{2} & g_{2} \vspace{1mm} \\
-e_{3} & -f_{3} & g_{3} & g_{3} & g_{3} \vspace{1mm} \\
-e_{3} & -f_{3} & g_{3} & g_{3} & g_{3} \vspace{1mm} \\
-e_{3} & -f_{3} & g_{3} & g_{3} & g_{3}
\end{array}
\right\|,
\end{equation}}
where
\begin{equation}\label{nopyt}
\begin{array}{lllll}
& a_{1}  = \frac{46\sqrt{34}}{561} - \frac{61}{132}, & a_{2} = - \frac{335\sqrt{34}}{4488} - \frac{5}{132}, & a_{3} = - \frac{20\sqrt{34}}{561} + \frac{5}{132},
\vspace{1mm} \\
& b_{1}  =  - \frac{29\sqrt{34}}{1122} + \frac{4}{33}, & b_{2} = \frac{95\sqrt{34}}{1122} + \frac{25}{66}, & b_{3} = \frac{37\sqrt{34}}{1122} + \frac{4}{33},
\vspace{1mm} \\
& c_{1} = \frac{7\sqrt{34}}{374} - \frac{5}{44}, & c_{2} = \frac{5\sqrt{34}}{1496} + \frac{5}{44}, & c_{3} = - \frac{\sqrt{34}}{1122} + \frac{7}{132},
\vspace{1mm} \\
& e_{1}  =\frac{46\sqrt{34}}{561} + \frac{61}{132}, & e_{2}  =   \frac{67\sqrt{34}}{4488} - \frac{1}{132}, & e_{3}  =\frac{20\sqrt{34}}{561} + \frac{5}{132},
\vspace{1mm} \\
& f_{1}  =-\frac{29\sqrt{34}}{1122} - \frac{4}{33}, & f_{2}  =\frac{95\sqrt{34}}{1122} - \frac{25}{66}, & f_{3}  = -\frac{37\sqrt{34}}{1122} + \frac{4}{33},
\vspace{1mm} \\
& g_{1}  = \frac{7\sqrt{34}}{374} + \frac{5}{44}, & g_{2}  = \frac{5\sqrt{34}}{1496} - \frac{5}{44}, & g_{3}  =\frac{\sqrt{34}}{1122} + \frac{7}{132}.
\end{array}
\end{equation}

Respectively, the matrix eigenvalues representation based on relation (\ref{decomba}) takes the following form, for $n \geq 1$,
\begin{align}\label{repres}
 \mathbf{P}_0^n & = \mathbf{\Pi} + (- \frac{1}{3})^n  \bar{\mathbf{\Pi}}_{2} +
(- \frac{1}{15} - \frac{1}{30} \sqrt{34})^n  \bar{\mathbf{\Pi}}_{3} \vspace{2mm} \nonumber \\
& \quad + (- \frac{1}{15} + \frac{1}{30} \sqrt{34})^n  \bar{\mathbf{\Pi}}_{4}.
\end{align}

Taking into account the specific forms of the damping matrix $\mathbf{D}$ and matrices  $\bar{\mathbf{\Pi}}_{2} = 2\mathbf{\Pi}_{2}$, we get that, in this case, coefficients  $\pi_{0, \bar{d}, j}[2] = 0, j = 1, \ldots, m$, and,
thus, the coefficients of the eigenvalues representation given in relation (\ref{decoma}) take the  following forms,  for $n \geq 1$,
\begin{equation}\label{coeff}
  p_{0,\bar{d}, j}(n)  = \left\{
\begin{array}{llll}
 \frac{5}{66}+\left(-\frac{223\sqrt{34}}{22440}+\frac{41}{660}\right) \left(- \frac{1}{15} - \frac{\sqrt{34}}{30}\right)^n \vspace{2mm} &  \\
 \quad + \left(-\frac{223\sqrt{34}}{22440}-\frac{41}{660}\right) \left(- \frac{1}{15} + \frac{\sqrt{34}}{30}\right)^n,     & \text{for} \ j = 1, \vspace{2mm}  \\
 \frac{8}{33}-\left(\frac{13\sqrt{34}}{5610}-\frac{7}{330}\right) \left(- \frac{1}{15} - \frac{\sqrt{34}}{30}\right)^n \vspace{2mm} &  \\
\quad + \left(\frac{13\sqrt{34}}{5610}+\frac{7}{330}\right)\left(- \frac{1}{15} + \frac{\sqrt{34}}{30}\right)^n,      & \text{for} \  j = 2, \vspace{2mm} \\
  \frac{5}{22}+\left(\frac{211\sqrt{34}}{36720}-\frac{37}{1080}\right)\left(- \frac{1}{15} - \frac{\sqrt{34}}{30}\right)^n \vspace{2mm} \\
\quad  + \left(\frac{19\sqrt{34}}{7480}+\frac{3}{220}\right) \left(- \frac{1}{15} + \frac{\sqrt{34}}{30} \right)^n ,      & \text{for} \   j = 3,\\
  \frac{5}{22}+\left(\frac{211\sqrt{34}}{36720}-\frac{37}{1080}\right)\left(- \frac{1}{15} - \frac{\sqrt{34}}{30}\right)^n \vspace{2mm} \\
\quad  + \left(\frac{19\sqrt{34}}{7480}+\frac{3}{220}\right) \left(- \frac{1}{15} + \frac{\sqrt{34}}{30} \right)^n ,      & \text{for} \   j = 4, \\
  \frac{5}{22}+\left(\frac{211\sqrt{34}}{36720}-\frac{37}{1080}\right)\left(- \frac{1}{15} - \frac{\sqrt{34}}{30}\right)^n \vspace{2mm} \\
\quad  + \left(\frac{19\sqrt{34}}{7480}+\frac{3}{220}\right) \left(- \frac{1}{15} + \frac{\sqrt{34}}{30} \right)^n ,      & \text{for} \   j = 5.
\end{array}
\right.
\end{equation}

Finally, the coefficients in the asymptotic expansion (\ref{srewavabe}) given by relation (\ref{last}) take the following form,  for all $j\in \XX$ and $n \geq 1$,

\begin{equation}\label{finasak}
\tilde{\pi}_{0,\bar{d}, 1}[n] = \left\{
\begin{array}{llll}
\frac{307}{2178}, \quad  & \text{for} \ n = 1, \vspace{2mm} \\
\left(\frac{307}{4356}-\frac{853}{74052}\sqrt{34}\right) \left(\frac{1}{33}+\frac{\sqrt{34}}{33}\right)^{n-1} \vspace{2mm} &  \\
 \quad +\left(\frac{307}{4356}+\frac{853}{74052}\sqrt{34}\right) \left(\frac{1}{33}-\frac{\sqrt{34}}{33}\right)^{n-1}, & \text{for} \ n > 1.
\end{array}
\right.
\end{equation}
\ \\
\begin{equation}\label{finasaka}
\tilde{\pi}_{0,\bar{d}, 2}[n] = \left\{
\begin{array}{llll}
 -\frac{50}{1089}, \quad  & \text{for} \ n = 1, \vspace{2mm}  \\
\left(-\frac{25}{1089}+\frac{107}{37026}\sqrt{34}\right) \left(\frac{1}{33}+\frac{\sqrt{34}}{33}\right)^{n-1} \vspace{2mm} &  \\
 \quad -\left(\frac{25}{1089}+\frac{107}{37026}\sqrt{34}\right) \left(\frac{1}{33}-\frac{\sqrt{34}}{33}\right)^{n-1}, & \text{for} \ n > 1.
\end{array}
\right.
\end{equation}
\ \\
\begin{align}\label{finasakak}
\tilde{\pi}_{0,\bar{d}, 3}[n]& = \tilde{\pi}_{0,\bar{d}, 4}[n] = \tilde{\pi}_{0,\bar{d}, 5}[n] \nonumber\\
 &=\left\{
\begin{array}{llll}
-\frac{23}{726}, \quad  & \text{for} \ n = 1, \vspace{2mm}  \\
\left(-\frac{23}{1452}+\frac{71}{24684}\sqrt{34}\right) \left(\frac{1}{33}+\frac{\sqrt{34}}{33}\right)^{n-1} \vspace{2mm} &  \\
 \quad -\left(\frac{23}{1452}+\frac{71}{24684}\sqrt{34}\right) \left(\frac{1}{33}-\frac{\sqrt{34}}{33}\right)^{n-1}, & \text{for} \ n > 1.
\end{array}
\right.
\end{align}

For example,  the second order expansions (\ref{srewavabe}), for stationary probabilities $\pi_{\e, j}, j = 1, \ldots, 5$, take the forms given below in relation (\ref{finasaba}). The terms of the expansions, except the stationary probabilities (first terms) which are exact values, are computed correct to 5 decimal digits.
\begin{equation}\label{finasaba}
\pi_{\e, j} \approx \left\{
\begin{array}{llll}
\frac{5}{66} + 0.14096 \e - 0.01946 \e^2 + O(\e^3) & \text{for} \ j = 1, \vspace{2mm}  \\
\frac{8}{33} - 0.04591 \e + 0.00456 \e^2 + O(\e^3) & \text{for} \ j = 2,  \vspace{2mm}  \\
\frac{5}{22} - 0.03168 \e + 0.00497 \e^2 + O(\e^3) & \text{for} \ j = 3,  \vspace{2mm} \\
\frac{5}{22} - 0.03168 \e + 0.00497 \e^2 + O(\e^3) & \text{for} \ j = 4,  \vspace{2mm} \\
\frac{5}{22} - 0.03168 \e + 0.00497 \e^2 + O(\e^3) & \text{for} \ j = 5.
\end{array}
\right.
\end{equation}




It can be observed from relation (\ref{finasaba}) that contributions of second terms in the above expansions are insignificant and, therefore, can be neglected even for moderate value of $\e = 0.2$. For example, the first term $0.14096 \e$ in the above asymptotic expansion for the stationary probability $\pi_{0, \bar{d}, 1} = \frac{5}{66} \approx 0.07576$ takes value $0.02820$ that is about $37.22\%$
of the corresponding  stationary probability. While, the  second term $0.01946 \e^2$ in the above asymptotic expansion takes value $0.00078$ that is only about $2.76\%$ of the first term and about $1.03\%$  of the corresponding  stationary probability. 

Let us also mention that the equalities  $\sum_{j \in \XX} \tilde{\pi}_{0, \bar{d}, j}[n] = 0, n \geq 1$
hold for coefficients given in relations (\ref{finasak})-(\ref{finasakak}) as well as, approximately up to the corresponding rounding corrections, for coefficients given in relation  (\ref{finasaba}).

Relation (\ref{decomba}) imply that inequality (\ref{asnao}) holds with the quantities,
$$
\lambda =  |\rho_{2}| \ {\rm and} \ C = \max_{r, k \in \XX} \sum_{l = 2}^{\bar{m}} | \bar{\pi}_{0, rk}[l] | \, | \frac{\bar{\rho}_{l}}{\bar{\rho}_{2}}|.
$$

The inequality (\ref{esti}) takes the following form,  for  $\e \in (0, 1]$ and $j  = 1, \ldots, 5$, 
\begin{equation}\label{estimabas}
| \pi_{\e, j}  - \pi_{0, j} | \leq \e \left(|\frac{1}{5} - \pi_{0, j}| +  \frac{67}{4488}\sqrt{34}+\frac{49}{132}\right), 
\end{equation}
where the stationary probabilities, $\pi_{0,  1} = \frac{5}{66}, \pi_{0, 2} = \frac{8}{33}, \pi_{0,  r} = \frac{5}{22}, r = 3, 4, 5$. 

Table \ref{Tab:bounds} give values $\pi_{\e,j}$, $\pi_{0,j}$, and $|\pi_{\e,j}-\pi_{0,j}|$ computed with the use of 
rounded solutions of the system of linear equations  (\ref{ergom}) and (\ref{ergsys}) and the upper bounds 
given by inequality (\ref{estimabas}), for the moderate value $\e = 0.15$,

 \begin{table}[H]
 \begin{center}
\begin{tabular}{c|c|c|c|c}
\hline 
$j \in \XX$ & $\pi_{\e,j} \vspace{1mm}$ & $\pi_{0,j}$ & $|\pi_{\e,j}-\pi_{0,j}|$ & Upper bound \\  \hline
1       & 0.09646 & 0.07576  & 0.02071 & 0.08738 \vspace{1mm}  \\
2       & 0.23564  &0.24242  & 0.00678 &  0.07510  \vspace{1mm}   \\
3,4,5 & 0.22263  & 0.22727 & 0.00464  & 0.07283   \vspace{1mm}  \\
\hline
\end{tabular}
\caption {Upper bounds for the rate of convergence  ($\e=0.15$).} \label{Tab:bounds}
\end{center}
\end{table}
\vspace{-4mm}

Figure \ref{finasaba1} displays  the asymptotical behaviour of the ergodicity coefficient $\Delta_N(\mathbf{P}_0)$ (computed with the use of relation (\ref{formades}),  as $N \to \infty$. In this case, coefficients of ergodicity $\Delta_N(\mathbf{P}_0) \to |\rho_2| = \frac{1}{3}$ as $N \to \infty$ consistently with relation (\ref{formadesam}). 

\begin{figure}[H]
  \centering
  \includegraphics[width=10cm]{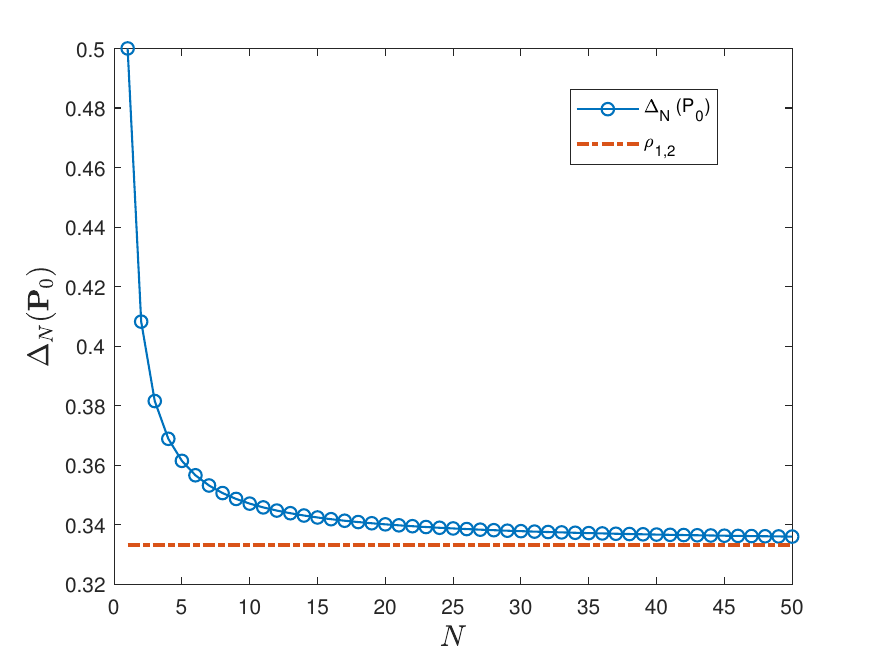}
  \caption{Ergodicity coefficient.}\label{finasaba1}
\end{figure}


\vspace{1mm}

{\bf 9.2 High order terms in the expansions for stationary probabilities of perturbed Markov chain}. Let us now present an the example, where the second term  in the corresponding asymptotic expansion for stationary probability can not  be ignored.

Here,  we consider a regularly perturbed Markov chain with matrices $\mathbf{P}_0$ and $\mathbf{D}$ given in relation (\ref{ex2}) and the link graph shown in Figure \ref{fig:example5}.  

 \begin{figure}[H]
\begin{center}
\includegraphics [width=.2\textwidth]{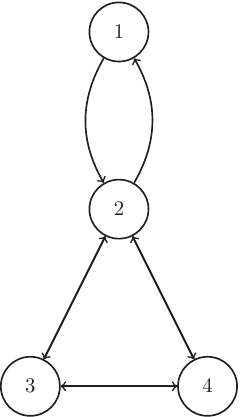}
\end{center}
\caption{A simple information network} \label{fig:example5}
\end{figure}

\begin{equation}\label{ex2}
\mathbf{P}_0 = \left\|
	\begin{array}{cccc}
0 & 1 & 0 & 0  \\
\frac{1}{3} & 0 & \frac{1}{3} & \frac{1}{3} \\
0 & \frac{1}{2} & 0 & \frac{1}{2}  \\
0 & \frac{1}{2} & \frac{1}{2} & 0
\end{array}
	\right\|,  \, \mathbf{D}   = \left\|
\begin{array}{cccc}
\frac{1}{4} & \frac{1}{4} & \frac{1}{4} & \frac{1}{4}  \vspace{1mm} \\
\frac{1}{4} & \frac{1}{4} & \frac{1}{4} & \frac{1}{4}  \vspace{1mm} \\
\frac{1}{4} & \frac{1}{4} & \frac{1}{4} & \frac{1}{4}  \vspace{1mm} \\
\frac{1}{4} & \frac{1}{4} & \frac{1}{4} & \frac{1}{4}
\end{array}
\right\|.
\end{equation}

\vspace{1mm}

In this case, condition ${\bf A}_{1}$ also holds.

The eigenvalues $\rho_k = \rho_{1, k} , k = 1, \ldots, 4$  of the matrix $\mathbf{P}_{0}$ are, 
\[\rho_{1}=1,\, \rho_{2}=-\frac{1}{4}-\frac{\sqrt{33}}{12},\,\rho_{3}=-\frac{1}{2}, \,  \rho_{4}=-\frac{1}{4}+\frac{\sqrt{33}}{12}.\]

Computations, analogous to performed above in the first example, give the following values for the coefficients of the eigenvalues 
representation given in relation (\ref{decomane}), for $n \geq 1$, 
\begin{equation}
  p_{0,\bar{d}, j}(n)  = \left\{
\begin{array}{llll}
 \frac{1}{8}+\left(\frac{1}{16}+\frac{\sqrt{33}}{528}\right)\left(-\frac{1}{4}-\frac{\sqrt{33}}{12}\right)^n \vspace{2mm} &  \\
 \quad + \left(\frac{1}{16}-\frac{\sqrt{33}}{528}\right)\left(-\frac{1}{4}+\frac{\sqrt{33}}{12}\right)^n,     & \text{for} \ j = 1, \vspace{2mm}  \\
 \frac{3}{8}-\left(\frac{1}{16}+\frac{3\sqrt{33}}{176}\right)\left(-\frac{1}{4}-\frac{\sqrt{33}}{12}\right)^n \vspace{2mm} &  \\
 \quad + \left(-\frac{1}{16}+\frac{3\sqrt{33}}{176}\right)\left(-\frac{1}{4}+\frac{\sqrt{33}}{12}\right)^n, & \text{for} \  j = 2, \vspace{2mm} \\
  \frac{1}{4}+\frac{\sqrt{33}}{132}\left(-\frac{1}{4}-\frac{\sqrt{33}}{12}\right)^n \vspace{2mm} \\
\quad  - \frac{\sqrt{33}}{132} \left(-\frac{1}{4}+\frac{\sqrt{33}}{12}\right)^n ,      & \text{for} \   j = 3,4.
\end{array}
\right.
\end{equation}

Finally, the coefficients in the asymptotic expansion (\ref{srewavabena}) given by relation (\ref{lastnem}) 
take the following form,  for $j \in \XX$,
\begin{equation}\label{finasan}
\tilde{\pi}_{0,\bar{d}, 1}[n] = \left\{
\begin{array}{llll}
\frac{7}{64}, \quad  & \text{for} \ n = 1, \vspace{2mm} \makebox[12mm]{} \\
\frac{(-1)^{2n+1}}{16}\left(\frac{3+\sqrt{33}}{15+\sqrt{33}}\right)^n \vspace{2mm} &  \\
 \quad  - \left(\frac{11+5\sqrt{33}}{176}\right)\left(\frac{-3+\sqrt{33}}{-15+\sqrt{33}}\right)^n \vspace{2mm} &  \\
 \quad +(-1)^{2n}\frac{5\sqrt{33}}{176} \left(\frac{3+\sqrt{33}}{15+\sqrt{33}}\right)^n, & \text{for} \ n > 1.
\end{array}
\right.
\end{equation}
\begin{equation}\label{finasana}
\makebox[9mm]{}\tilde{\pi}_{0,\bar{d}, 2}[n] = \left\{
\begin{array}{llll}
 - \frac{3}{64}, \quad  & \text{for} \ n = 1, \vspace{1mm} \\
(-1)^{1+2n}\left(\frac{33+\sqrt{33}}{176}\right)\left(\frac{3+\sqrt{33}}{15+\sqrt{33}}\right)^n \vspace{2mm} \makebox[12mm]{} &  \\
 \quad +\left(\frac{-33+\sqrt{33}}{176}\right) \left(\frac{-3+\sqrt{33}}{-15+\sqrt{33}}\right)^n, & \text{for} \ n > 1.
\end{array}
\right.
\end{equation}

\begin{align}\label{finasanan}
\tilde{\pi}_{0,\bar{d}, 3}[n]& = \tilde{\pi}_{0,\bar{d}, 4}[n] \nonumber\\
 &=\left\{
\begin{array}{llll}
 - \frac{1}{32}, \quad  & \text{for} \ n = 1, \vspace{1mm} \\
\frac{\sqrt{33}}{11(3+\sqrt{33})}(-1)^{n}\left(\frac{-3-\sqrt{33}}{15+\sqrt{33}}\right)^n \vspace{2mm} &  \\
 \quad +\frac{\sqrt{33}}{11(-3+\sqrt{33})}(-1)^{n}\left(\frac{-3+\sqrt{33}}{15-\sqrt{33}}\right)^n, & \text{for} \ n > 1. 
\end{array}
\right.
\end{align}

Consequently, the second order asymptotic expansions (\ref{srewavabena}), for stationary probabilities $\pi_{\e, j}, j = 1, \ldots, 4$, take the forms given below in relation (\ref{expa2}). The terms of the expansions, except the stationary probabilities (first terms) which are exact values, are computed correct to 5 decimal digits.

\begin{equation}\label{finasanana}
\makebox[11mm]{}   \pi_{\e, j} = \left\{
 \begin{array}{llll}
\frac{1}{8}+\frac{7}{64}\e-\frac{1}{512}\e^{2} +O(\e^{3}), \quad  & \text{for} \ j = 1, \makebox[20mm]{} \vspace{1mm}\\
\frac{3}{8}-\frac{3}{64}\e-\frac{27}{512}\e^{2} +O(\e^{3}), \quad  & \text{for} \ j = 2, \vspace{1mm}\\
\frac{1}{4}-\frac{1}{32}\e+\frac{7}{256}\e^{2} +O(\e^{3}),  \quad  & \text{for} \ j = 3, \vspace{1mm}\\
\frac{1}{4}-\frac{1}{32}\e+\frac{7}{256}\e^{2} +O(\e^{3}),  \quad  & \text{for} \ j = 4.
\end{array}
\right.
 \end{equation}
\begin{equation}\label{expa2}
\makebox[14mm]{} \approx \left\{
  \begin{array}{llll}
    \frac{1}{8}+ 0.10937 \e- 0.00195 \e^{2} +O(\e^{3}), \quad  & \text{for} \ j = 1,  \vspace{1mm}\\
    \frac{3}{8}- 0.04687 \e- 0.05273 \e^{2} +O(\e^{3}), \quad  & \text{for} \ j = 2, \vspace{1mm}\\
    \frac{1}{4}- 0.03125 \e+ 0.02734 \e^{2} +O(\e^{3}),  \quad  & \text{for} \ j = 3, \vspace{1mm}\\
    \frac{1}{4}- 0.03125\e+ 0.02734 \e^{2} +O(\e^{3}),  \quad  & \text{for} \ j = 4.
  \end{array}
  \right.
\end{equation}

Let us mention again that equalities $\sum_{j \in \XX} \tilde{\pi}_{0, \bar{d}, j}[n] = 0, n \geq 1$
 hold for coefficients given in relations 
(\ref{finasan})-(\ref{finasanan}) and (\ref{finasanana}) as well as, approximately up to the corresponding rounding corrections, for coefficients given in  relation  (\ref{expa2}).

We see that in the third ($j = 3$) and fourth ($j = 4$) expansions in relation (\ref{expa2})  for $\e = 0.2$, the term {\footnotesize$\frac{1}{32}\e$} is about 2.5\% of the limiting stationary probability 0.25, while the term {\footnotesize $\frac{7}{256} \e^2$} is about 0.4\% of 0.25, i.e. it improves the accuracy of approximation for the stationary probability for about 16\%. \\

\vspace{1mm}

{\bf 9.3 Singularly perturbed Markov chains with damping component in the triangular array mode}.
In this example we consider an ergodic singularly perturbed Markov chains  with damping component as stated in Theorem 9. We use a two-disjoint information networks given in Figure \ref{fig:example3}, whose matrices $\mathbf{P}_0$ and $\mathbf{D}$ are given by (\ref{eqn.ex3}).

\begin{figure}[H]
\begin{center}
\includegraphics [width=.4\textwidth]{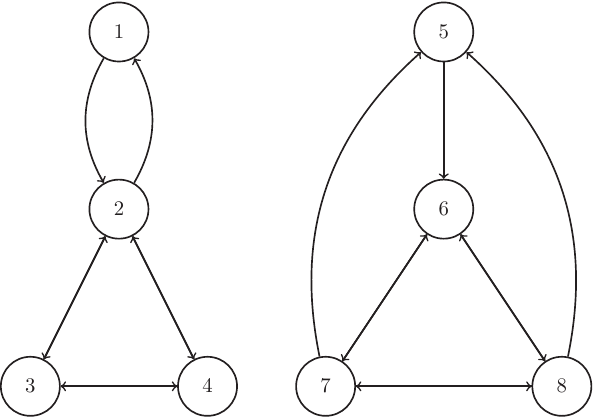}
\end{center}
\caption{A simple two disjoint information networks} \label{fig:example3}
\end{figure}
{\footnotesize
\begin{equation}\label{eqn.ex3}
\mathbf{P}_{0}=\left\|
	\begin{array}{cccccccc}
0 & 1 & 0 & 0 & 0 & 0 & 0 & 0 \\
\frac{1}{3} & 0 & \frac{1}{3} & \frac{1}{3} & 0 & 0 & 0 &0 \vspace{1mm} \\
0 & \frac{1}{2} & 0 & \frac{1}{2} & 0 & 0 & 0 & 0 \vspace{1mm} \\
0 & \frac{1}{2} & \frac{1}{2} & 0 & 0 & 0& 0 & 0 \vspace{1mm} \\
0 & 0 & 0 & 0& 0 & 1 & 0 & 0 \vspace{1mm} \\
0 & 0 & 0 & 0&0& 0  & \frac{1}{2}  & \frac{1}{2}  \vspace{1mm}  \\
0 & 0 & 0 & 0 & \frac{1}{3}  & \frac{1}{3}  & 0& \frac{1}{3}  \vspace{1mm} \\
0 & 0 & 0 & 0 & \frac{1}{3}  & \frac{1}{3}  & \frac{1}{3} & 0

\end{array}
	\right\|,  \, \mathbf{D}   = \left\|
	\begin{array}{cccccccc}
\frac{1}{8} & \frac{1}{8} & \frac{1}{8} & \frac{1}{8} & \frac{1}{8} & \frac{1}{8} & \frac{1}{8} & \frac{1}{8} \vspace{1mm} \\
\frac{1}{8} & \frac{1}{8} & \frac{1}{8} & \frac{1}{8} & \frac{1}{8} & \frac{1}{8} & \frac{1}{8} & \frac{1}{8}  \vspace{1mm} \\
\frac{1}{8} & \frac{1}{8} & \frac{1}{8} & \frac{1}{8} & \frac{1}{8} & \frac{1}{8} & \frac{1}{8} & \frac{1}{8}  \vspace{1mm} \\
\frac{1}{8} & \frac{1}{8} & \frac{1}{8} & \frac{1}{8} & \frac{1}{8} & \frac{1}{8} & \frac{1}{8} & \frac{1}{8}  \vspace{1mm} \\
\frac{1}{8} & \frac{1}{8} & \frac{1}{8} & \frac{1}{8} & \frac{1}{8} & \frac{1}{8} & \frac{1}{8} & \frac{1}{8}  \vspace{1mm} \\
\frac{1}{8} & \frac{1}{8} & \frac{1}{8} & \frac{1}{8} & \frac{1}{8} & \frac{1}{8} & \frac{1}{8} & \frac{1}{8}  \vspace{1mm} \\
\frac{1}{8} & \frac{1}{8} & \frac{1}{8} & \frac{1}{8} & \frac{1}{8} & \frac{1}{8} & \frac{1}{8} & \frac{1}{8}  \vspace{1mm} \\
\frac{1}{8} & \frac{1}{8} & \frac{1}{8} & \frac{1}{8} & \frac{1}{8} & \frac{1}{8} & \frac{1}{8} & \frac{1}{8}

\end{array}
	\right\|
\end{equation}}

Here we define a phase space $\XX$ corresponding to $\mathbf{P}_0$ as $\XX =\XX^{(1)} \cup \XX^{(2)}$, where $\XX^{(1)}=\{1,2,3,4\}$ and $\XX^{(2)}=\{5,6,7,8\}$. The corresponding matrices of transition probabilities for the two ergodic classes are, respectively, $\mathbf{P}_{0,1}$ and $\mathbf{P}_{0,2}$, where
\[
   \mathbf{P}_{0,1}= \left\|
	\begin{array}{cccc}
	 0 & 1 & 0 & 0 \vspace{1mm} \\
	 \frac{1}{3} & 0 & \frac{1}{3} & \frac{1}{3} \vspace{1mm} \\
	 0 & \frac{1}{2} & 0 & \frac{1}{2} \vspace{1mm} \\
	 0 & \frac{1}{2} & \frac{1}{2} & 0
	\end{array}
	\right\|,  \, \mathbf{P}_{0,2}= \left\|
	\begin{array}{cccc}
	 0 & 1 & 0 & 0 \vspace{1mm} \\
     0 & 0  & \frac{1}{2}  & \frac{1}{2}  \vspace{1mm}  \\
    \frac{1}{3}  & \frac{1}{3}  & 0& \frac{1}{3}  \vspace{1mm} \\
    \frac{1}{3}  & \frac{1}{3}  & \frac{1}{3} & 0
	\end{array}
	\right\|.
\]

In this case,  condition ${\bf B}_{1}$ holds, with  $h = 2$.

We are going to get the second order asymptotic expansion for stationary probabilities  $\pi_{\e,k}, k\in \XX$ given by relation (\ref{srewavabena}). In fact, we can use Theorem 4 and to give the corresponding expansions for stationary probabilities  $\pi_{\e,k}^{(1)}, k \in \XX^{(1)}$ and $\pi_{\e,k}^{(2)}, k \in \XX^{(2)}$. Then, these expansions can be transformed  in the corresponding expansions for $\pi_{\e,k}, k\in \XX$, using relations $\pi_{\e,k} = 
f_{\bar{d}}^{(g)}\pi_{\e,k}^{(g)}, k \in \XX^{(g)}, g = 1, 2$ and taking into account that, in this case,
probabilities $f_{\bar{d}}^{(g)}=\sum_{k \in \XX^{(g)}}d_k = \frac{1}{2}, g = 1, 2$.

One can easily notice that matrix $\mathbf{P}_{0,1}$ coincides with matrix $\mathbf{P}_0$ given in relation (\ref{ex2}). Thus, the second order asymptotic expansions for stationary probabilities $\pi_{\e,k}^{(1)}, k \in \XX^{(1)}$ take forms given in relation (\ref{expa2}), and, thus,  the corresponding asymptotic expansions for stationary probabilities $\pi_{\e,k} = \frac{1}{2}\pi_{\e,k}^{(1)}, k \in \XX^{(1)}$ take the following forms, 
\begin{equation}\label{expa3a}
   \pi_{\e,k} =
  \begin{cases}
    \frac{1}{16}+\frac{7}{128}\e-\frac{1}{1024}\e^{2} +O(\e^{3}), \quad  & \text{for} \ k=1, \vspace{1mm} \\
    \frac{3}{16}-\frac{3}{128}\e-\frac{27}{1024}\e^{2} +O(\e^{3}), \quad  & \text{for} \ k=2, \vspace{1mm} \\
    \frac{1}{8}-\frac{1}{64}\e+\frac{7}{512}\e^{2} +O(\e^{3}), \quad  & \text{for} \ k=3, \vspace{1mm}  \\
    \frac{1}{8}-\frac{1}{64}\e+\frac{7}{512}\e^{2} +O(\e^{3}), \quad  & \text{for} \ k=4.
  \end{cases}
\end{equation}

 Similarly, the asymptotic expansion $\pi_{\e,k}^{(2)}$ for phase space $\XX^{(2)}$ can be computed. The eigenvalues corresponding to $\mathbf{P}_{0,2}$ are,
 \[\rho_{2, 1}=1,\, \rho_{2, 2}=-\frac{1}{3}-\frac{\sqrt{2}}{3}\mathbf{i},\,\rho_{2, 3}=-\frac{1}{3}+\frac{\sqrt{2}}{3}\mathbf{i}, \, \rho_{2, 4}=-\frac{1}{3},\]
 where $\mathbf{i}^{2} = -1$. Hence, from relation (\ref{decomane}), we obtain eigenvalues decomposition of $\mathbf{P}_{0,2}$, for $k \in \XX^{(2)}, n \geq 1$, as
\begin{equation}
  p_{0,\bar{d}, k}(n)  = \left\{
\begin{array}{llll}
 \frac{1}{6}+\left(\frac{\sqrt{2}\mathbf{i}}{16(1+\sqrt{2}\mathbf{i})}\right)\left(-\frac{1}{3}-\frac{\sqrt{2}}{3}\mathbf{i}\right)^n \vspace{2mm} &  \\
 \quad + \left(\frac{(-1+2\sqrt{2}\mathbf{i})}{24(-2+\sqrt{2}\mathbf{i})}\right)\left(-\frac{1}{3}+\frac{\sqrt{2}}{3}\mathbf{i}\right)^n,     & \text{for} \ k = 5, \vspace{2mm}  \\
 \frac{1}{3}-\left(\frac{1}{8+\sqrt{128}\mathbf{i}}\right)\left(-\frac{1}{3}-\frac{\sqrt{2}}{3}\mathbf{i}\right)^n \vspace{2mm} &  \\
 \quad - \left(\frac{-4-\sqrt{2}\mathbf{i}}{24(-2+\sqrt{2}\mathbf{i})}\right)\left(-\frac{1}{3}+\frac{\sqrt{2}}{3}\mathbf{i}\right)^n, & \text{for} \  k = 6, \vspace{2mm} \\
  \frac{1}{4}-\frac{\sqrt{2}\mathbf{i}}{32}\left(-\frac{1}{3}-\frac{\sqrt{2}}{3}\mathbf{i}\right)^n \vspace{2mm} \\
\quad  + \frac{\sqrt{2}\mathbf{i}}{32} \left(-\frac{1}{3}+\frac{\sqrt{2}}{3}\mathbf{i}\right)^n ,      & \text{for} \   k = 7, 8.
\end{array}
\right.
\end{equation}

Next, the coefficients $\tilde{\pi}^{(2)}_{0,\bar{d}, k}[n], k \in \XX^{(2)}$ (computed from relation (\ref{lastnem})), for $n\geq1$, take the form,
\begin{equation}\label{finasal}
\tilde{\pi}^{(2)}_{0,\bar{d}, 5}[n] = \left\{
\begin{array}{llll}
\frac{5}{72}, \quad  & \text{for} \ n = 1, \vspace{1mm} \\
\left(\frac{5}{144}+\frac{\sqrt{2}\mathbf{i}}{144}\right) \left(\frac{1}{3}+\frac{\sqrt{2}\mathbf{i}}{6}\right)^{n-1} \vspace{2mm} &  \\
 \quad +\left(\frac{5}{144}-\frac{\sqrt{2}\mathbf{i}}{144}\right) \left(\frac{1}{3}-\frac{\sqrt{2}\mathbf{i}}{6}\right)^{n-1}, & \text{for} \ n > 1.
\end{array}
\right.
\end{equation}
\begin{equation}\label{finasala}
\tilde{\pi}^{(2)}_{0,\bar{d}, 6}[n] = \left\{
\begin{array}{llll}
 - \frac{1}{36}, \quad  & \text{for} \ n = 1, \vspace{1mm} \\
\left(-\frac{1}{72}+\frac{5\sqrt{2}\mathbf{i}}{144}\right) \left(\frac{1}{3}+\frac{\sqrt{2}\mathbf{i}}{6}\right)^{n-1} \vspace{2mm} &  \\
 \quad -\left(\frac{1}{72}+\frac{5\sqrt{2}\mathbf{i}}{144}\right) \left(\frac{1}{3}-\frac{\sqrt{2}\mathbf{i}}{6}\right)^{n-1}, & \text{for} \ n > 1.
\end{array}
\right.
\end{equation}
\begin{align}\label{finasalal}
\tilde{\pi}^{(2)}_{0,\bar{d}, 7}[n]&= \tilde{\pi}^{(2)}_{0,\bar{d}, 8}[n] \nonumber\\
 &=\left\{
\begin{array}{llll}
 - \frac{1}{48}, \quad  & \text{for} \ n = 1, \vspace{1mm}  \\
\left(-\frac{1}{96}-\frac{\sqrt{2}\mathbf{i}}{48}\right) \left(\frac{1}{3}+\frac{\sqrt{2}\mathbf{i}}{6}\right)^{n-1} \vspace{2mm} &  \\
 \quad +\left(-\frac{1}{96}+\frac{\sqrt{2}\mathbf{i}}{48}\right) \left(\frac{1}{3}-\frac{\sqrt{2}\mathbf{i}}{6}\right)^{n-1}, & \text{for} \ n > 1.
\end{array}
\right.
\end{align}

Note that, according to remarks made in Subsection 6.1, coefficients $\tilde{\pi}^{(2)}_{0,\bar{d}, k}[n]$ are real numbers. They are 
given in relations (\ref{finasal})-(\ref{finasalal}) in the compact form, which can create a false impression that these coefficients may be 
complex numbers. However, one can see that the use  Newton binomial formula and following multiplication in the corresponding expressions 
yields cancelation of all complex terms. 

The asymptotic expansions for  $\pi_{\e,k} = \frac{1}{2}\pi^{(2)}_{\e,k}, k\in \XX^{(2)}$ take the following forms, 
\begin{equation}\label{expa3bb}
   \pi_{\e,k} =\left\{
       \begin{array}{llll}
    \frac{1}{12} +\frac{5}{144}\e  +\frac{1}{108}\e^{2}  +O(\e^{3}),\quad  & \text{for} \ k = 5, \vspace{2mm} \\
    \frac{1}{6}  -  \frac{1}{72}\e  -\frac{7}{432}\e^{2}  +O(\e^{3}),\quad  & \text{for} \ k= 6, \vspace{2mm} \\
    \frac{1}{8}  -  \frac{1}{96}\e  +\frac{1}{288}\e^{2}  + O(\e^{3}),\quad  & \text{for} \ k= 7,\vspace{2mm}  \\
    \frac{1}{8}  -\frac{1}{96}\e  +\frac{1}{288}\e^{2}  + O(\e^{3}),\quad  & \text{for} \ k= 8.
  \end{array}
\right.
\end{equation}

Hence, by combining relations (\ref{expa3a}) and (\ref{expa3bb}), we obtain the asymptotic expansions for  stationary probabilities  $\pi_{\e,k}, k \in \XX = \XX^{(1)} \cup \XX^{(2)}$,
\begin{equation}\label{expa3allb}
   \pi_{\e,k} =\left\{
    \begin{array}{llll}
     \frac{1}{16}+\frac{7}{128}\e-\frac{1}{1024}\e^{2} +O(\e^{3}),\quad  & \text{for} \ k = 1, \vspace{2mm} \\
       \frac{3}{16}-\frac{3}{128}\e-\frac{27}{1024}\e^{2} +O(\e^{3}),\quad  & \text{for} \ k = 2, \vspace{2mm} \\
       \frac{1}{8}-\frac{1}{64}\e+\frac{7}{512}\e^{2} +O(\e^{3}),\quad  & \text{for} \ k = 3, \vspace{2mm} \\
       \frac{1}{8}-\frac{1}{64}\e+\frac{7}{512}\e^{2} +O(\e^{3}),\quad  & \text{for} \ k = 4. \vspace{2mm} \\
       \frac{1}{12} +\frac{5}{144}\e  +\frac{1}{108}\e^{2}  +O(\e^{3}),\quad  & \text{for} \ k = 5, \vspace{2mm} \\
       \frac{1}{6}  -  \frac{1}{72}\e  -\frac{7}{432}\e^{2}  +O(\e^{3}),\quad  & \text{for} \ k = 6, \vspace{2mm} \\
       \frac{1}{8}  -  \frac{1}{96}\e  +\frac{1}{288}\e^{2}  + O(\e^{3}),\quad  & \text{for} \ k = 7, \vspace{2mm} \\
       \frac{1}{8}  -\frac{1}{96}\e  +\frac{1}{288}\e^{2}  + O(\e^{3}),\quad  & \text{for} \ k = 8. \vspace{2mm}  \\
          \end{array}
\right.
\end{equation}

Let us mention ones more time  that equalities $\sum_{k \in \XX} \tilde{\pi}_{0, \bar{d}, k}[n] = 0, n \geq 1$  hold for coefficients given in relations (\ref{finasal})-(\ref{finasalal}) and (\ref{expa3allb}). 

Expansions (\ref{expa3allb}) let us illustrate the results  presented in Theorem 9. The most interesting in this theorem is the asymptotic relation (\ref{firsta}).






In order to illustrate this relation,  we choose the initial distribution $\bar{p}= \langle1,0,0,0,0,0,0,0 \rangle$, which  is concentrated in one ergodic class $\XX^{(1)}$. In this case, probabilities $f^{(1)}_{\bar{p}}=1$, $f^{(2)}_{\bar{p}} = 0$. Note that  probabilities $f^{(1)}_{\bar{d}} = f^{(2)}_{\bar{d}} = \frac{1}{2}$ for the initial distribution $\bar{p} = \bar{d}$.   Probabilities $\pi_{0,\bar{p},k}$  and $\pi_{0,\bar{d}, k}$  are computed using relation (\ref{convert}). 
The  asymptotic relation  (\ref{firsta}) takes, for $k = 1$, the form, $p_{\e,\bar{p}, 1}(n_{\e}) \to \pi_{0,\bar{p}, 1}(t)$ for $n_{\e}\to \infty$, $\e n_{\e} \to t$ and $\e \to 0$.  Further, we choose $\e = 0.1$. According to the asymptotic relation (\ref{firsta}), the values $p_{\e,\bar{p}, 1}(n_\e)$, for $n_\e$ such that $\e n_\e \approx t$ can be expected to take values close to $\pi_{0,\bar{p}, 1}(t)$, for $1 \leq t \leq 3$. This,  indeed, can be
seen in Figure \ref{finasaba2b}.

In addition to that result shows that the relative absolute errors, $|\pi_{0,\bar{p},1}(t)-\pi_{0,\bar{d},1}|$ decreases dramatically from about 37\% to 4\%  as $\e n_{\e}$  increases from 1 to 3 respectively.

\begin{figure}[H]
  \centering
  \includegraphics[width=12cm]{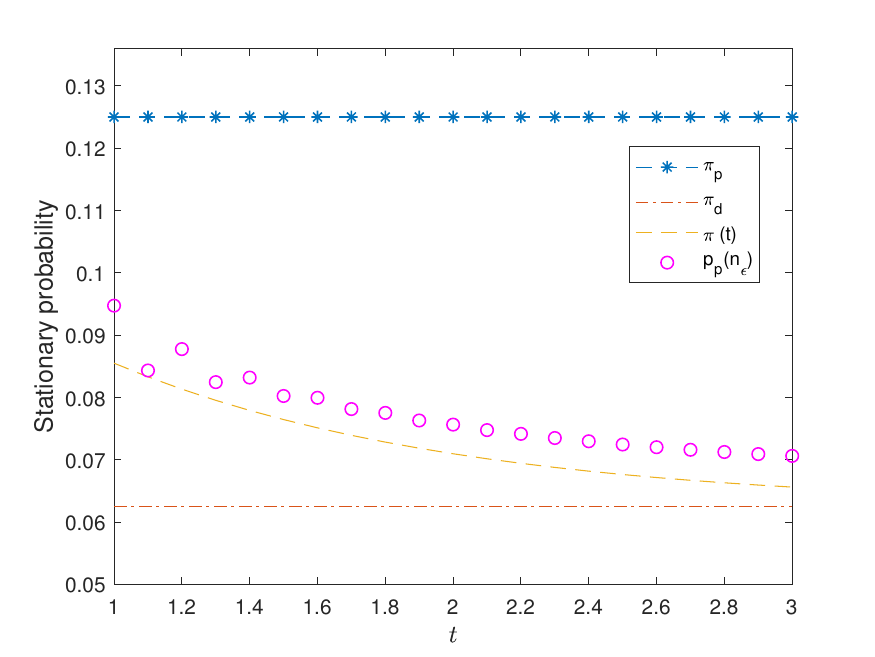}
  \caption{An illustration of Theorem 8 (ii), where $\pi_{p}$, $p_{p}(n_\e)$, $\pi (t)$ and $\pi_{d}$ represent $\pi_{0,\bar{p},1}$,  $p_{\e,\bar{p},1}(n_{\e})$, $\pi_{0,\bar{p},1}(t)$ and $p_{\e,\bar{d},1}$ respectively.} \label{finasaba2b}
\vspace{-3mm}
\end{figure}


{\bf 9.4  Concluding remarks and comments}. One of the main reasons for approximation of  the Markov chain $X_{0, n}$, with the matrix of transition probabilities
$\mathbf{P}_0$,  by perturbed (regularised) Markov chains with damping component $X_{\e, n}$, with the matrix of transition probabilities
$\mathbf{P}_\e =  (1 - \e)\mathbf{P}_0   + \e \mathbf{D}$, is to use it for approximation of the stationary distribution $\bar{\pi}_{0} = \langle \pi_{0, j}, j \in \XX \rangle$ of  the Markov chain $X_{0, n}$ by the stationary distribution $\bar{\pi}_{\e} = \langle \pi_{\e, j},
j \in \XX \rangle$ of  the Markov chain $X_{\e, n}$.

Since the corresponding phase space $\XX = \{ 1, \ldots, m \}$ can be large, the power method can be used for approximative computing of stationary
distribution $\bar{\pi}_{\e}$, where its components $\pi_{\e, j}$ are approximated by probabilities  $p_{\e, \bar{p}, j}(n) =
\sum_{i \in \XX} p_i p_{\e, ij}(n)$, where $p_{\e, ij}(n)$ are elements of the matrix $\mathbf{P}_\e^n$ and $\bar{p} =  \langle p_{j},
j \in \XX \rangle$ is some  initial distribution.

The results given in Theorems 1 -- 9 show that the situation significantly  differs for two models, (a) where the phase space $\XX$  is one class of communicative states for the Markov chain $X_{0, n}$, i.e., this Markov chain is ergodic,  and, (b) where the phase space $\XX$  splits in several closed classes of communicative states for the Markov chain $X_{0, n}$, and, thus, this Markov chain is not ergodic.

Rates of convergence of stationary probabilities $\pi_{\e, j}$ to $\pi_{0, j}$ as $\e \to 0$ and probabilities  $p_{\e, \bar{p}, j}(n)$ to stationary probabilities  $\pi_{0, \bar{p}, j}$ as $n \to \infty$ play the key role  in the above method.

We give  explicit upper bounds  for  rates of convergence
of stationary probabilities $\pi_{\e, j}$ to  $\pi_{0, j}$ and asymptotic expansions for  $\pi_{\e, j}$, with respect to damping parameter $\e$, in Theorems 1 -- 4.  We also give explicit upper bounds  for  rates of convergence of probabilities $p_{\e, \bar{p}, j}(n)$ to stationary probabilities $\pi_{\e, \bar{p}, j}$, in Theorems 5 -- 7.

These results let one construct the two-stage effective algorithms for approximating the stationary distribution $\bar{\pi}_{0}$ in the case,  where the phase space $\XX$  is one class of communicative states for the Markov chain $X_{0, n}$, and, thus, this Markov chain is ergodic.

At the first stage,  one approximates the stationary probabilities $\pi_{\e, j}$ by probabilities  $p_{\e, \bar{p}, j}(n)$.  The rate of this approximation has the order $O( \Delta_N(\mathbf{P}_0)^n (1 -  \e)^n)$. The effectiveness of this approximation declines for small values of damping parameter $\e$ and values of ergodicity coefficient $\Delta_N(\mathbf{P}_0)$ closed to $1$ that is typical for models with the phase space $\XX$  nearly decomposed in several closed
classes of states.

At the second stage, one approximate stationary probabilities  $\pi_{0, j}$ by stationary probabilities $\pi_{\e, j}$. The rate of approximation has the order   $O( \e)$, which also, can be improved by using the corresponding asymptotic expansions.

Here, some dual effect takes place. At the second stage it would be better to use small values of regularisation parameter $\e$, while at the first stage using small values of $\e$ are not desirable.

Theorem 8 let one also to use the one-step variant of the algorithm described above and approximate the stationary probabilities $\pi_{0, j}$ by probabilities $p_{\e, \bar{p}, j}(n_\e)$, with an arbitrary positive integers $n_\e \to \infty$ as $\e \to 0$. Moreover, the explicit upper bounds for $|p_{\e, \bar{p}, j}(n_\e) - \pi_{0, j}|$ pointed out in Remark 6 let one balance the choice of $\e$ and $n_\e$.

The case,  where the phase space $\XX$ splits in several closed classes of communicative states for the Markov chain $X_{0, n}$, and, thus, this Markov chain is not ergodic, is more complex.

As a matter of fact, in this case,  (c) stationary probabilities $\pi_{0, \bar{p}, k}$  for the Markov chain $X_{0, n}$ depend on the initial distribution $\bar{p}$,  (d) the stationary probabilities $\pi_{\e, k}$ for the Markov chain $X_{\e, n}$ converge to  the stationary probabilities  $\pi_{0,  \bar{d}, k}$ as $\e \to 0$.

If the initial distribution $\bar{p} = \bar{d}$ then the two-stage algorithm  as well as its one-stage variant described above can be applied.

However, if $\bar{p} \neq \bar{d}$, one should be more careful, since in this case  it may be that the stationary probability
$\pi_{0, \bar{d}, k} \neq \pi_{0, \bar{p}, k}$ and, thus, the two-stage algorithm described above does not work.

In this case, Theorem 9  answers the question about applicability quantities  $p_{\e, \bar{p}, j}(n_\e)$ as approximations for
stationary probabilities for the Markov chain $X_{0, n}$. In fact, these probabilities converge to some mixture
of stationary probabilities $\pi_{0, \bar{p}, k}$ and $\pi_{0, \bar{d}, k}$, namely,  $\pi_{0, \bar{p}, k}(t) = \pi_{0, \bar{p}, k} e^{- t}  + \pi_{0, \bar{d},  k}(1 - e^{- t})$,  as $n_\e \to \infty$ in such way that
$\e n_\e \to t \in [0, \infty]$ as $\e \to 0$. Moreover, the explicit upper bounds for $|p_{\e, \bar{p}, k}(n_\e) -
\pi_{0, \bar{p}, k}(t)|$ pointed out in Remark 7 let one balance
the choice of $\e$ and $n_\e$ and, in some sense, predict the value of limit $\pi_{0, \bar{p}, k}(t)$ depending on the value of
quantity $\e n_\e$.

The computational examples presented in Section 9 illustrate  in which way results given in Theorems 1 - 9 can be used in experimental  studies of Markov chains with damping components associated with information networks, in particular, in studies of PageRank algorithms.

In conclusion, we would like to note that Theorems 1 -- 9 present results of perturbation analysis for Markov chains with damping component for the basic cases, where the phase space of the unperturbed Markov chain $X_{0, n}$  either is one aperiodic class of communicative states or split in several closed aperiodic classes of communicative states. Despite some technical complications, analogous results  can be also obtained for the cases, where the phase space of the unperturbed Markov chain $X_{0, n}$  either is one periodic class of communicative states or split in several closed periodic classes of communicative states plus possibly a class of transient states.  We are going to present such results in future publications as well as results of the corresponding experimental  studies. \\

\noindent{\bf Acknowledgements} \\

This research was supported by the Swedish International Development Cooperation Agency (Sida), International Science Programme (ISP) in Mathematical Sciences
(IPMS) and Sida Bilateral Research Programmes for research and education capacity development in Mathematics in Uganda and Tanzania. The authors are also grateful to the research environment Mathematics and Applied Mathematics (MAM), Division of Applied Mathematics, M\"{a}lardalen University for providing an excellent and inspiring environment for research education and research.

\end{document}